\documentclass[a4paper,12pt]{article}
\usepackage{amssymb}
\usepackage{amsmath}
\usepackage{amsthm}
\usepackage{graphicx}
\numberwithin{equation}{section}

\newcommand{\sumn}{\sum_{n=1}^{\infty}}
\newcommand{\sumj}{\sum_{j=1}^{\infty}}
\newcommand{\ep}{\varepsilon}
\newcommand{\la}{\lambda}
\newcommand{\va}{\varphi}
\newcommand{\ppp}{\partial}

\newcommand{\www}{\widetilde}

\newcommand{\R}{\mathbb{R}}
\newcommand{\Q}{\mathbb{Q}}
\newcommand{\C}{\mathbb{C}}
\newcommand{\N}{\mathbb{N}}

\newcommand{\ooo}{\overline}
\newcommand{\OOO}{\Omega}

\allowdisplaybreaks
\title{Determination of source terms in diffusion and
wave equations by observations after incidents:
uniqueness and stability}
\author{
Jin Cheng\footnotemark[1], \quad Shuai Lu\footnotemark[1],  \quad
Masahiro Yamamoto \footnotemark[2]\, \footnotemark[3]\, \footnotemark[4]\, \footnotemark[5]
}
\footnotetext[1]{School of Mathematical Sciences, Fudan University, Shanghai, China, {\tt jcheng@fudan.edu.cn; slu@fudan.edu.cn}.}
\footnotetext[2]{Department of Mathematical Sciences, The University of Tokyo, Komaba, Meguro, Tokyo 153, Japan,
{\tt myama@ms.u-tokyo.ac.jp}}
\footnotetext[3]{Honorary Member of Academy of Romanian Scientists,
Ilfov, nr. 3, Bucuresti, Romania}
\footnotetext[4]{Correspondence member of Accademia Peloritana dei Pericolanti
}
\footnotetext[5]{Peoples' Friendship University of Russia
(RUDN University) 6 Miklukho-Maklaya St, Moscow, 117198, Russian Federation
}

\begin{document}
\maketitle
\newpage

\begin{abstract}
We consider a diffusion and a wave equations:
$$
\ppp_t^ku(x,t) = \Delta u(x,t) + \mu(t)f(x), \quad
x\in \OOO, \, t>0,  \quad k=1,2
$$
with the zero initial and boundary conditions, where
$\OOO \subset \R^d$ is a bounded domain.  We establish
uniqueness and/or stability results for inverse problems of
\begin{itemize}
\item
determining $\mu(t)$, $0<t<T$ with given $f(x)$.
\item
determining $f(x)$, $x\in \OOO$ with given $\mu(t)$
\end{itemize}
by data of $u$: $u(x_0,\cdot)$ with fixed point
$x_0\in \OOO$ or Neumann data on subboundary over time
interval.   In our inverse problems, data are taken over time
interval $T_1<t<T_1$,
by assuming that $T<T_1<T_2$ and $\mu(t)=0$ for $t\ge T$, which means that
the source stops to be active after the time $T$ and the observations
are started only after $T$.  This assumption is practical by
such a posteriori data after incidents,
although inverse problems had been well studied in the case of $T=0$.
We establish the non-uniqueness, the uniqueness and
conditional stability for a diffusion and a wave equations.
The proofs are based on eigenfunction expansions of the solutions $u(x,t)$,
and we rely on various knowledge of
the generalized Weierstrass theorem on polynomial approximation,
almost periodic functions, Carleman estimate, non-harmonic
Fourier series.
\end{abstract}

\baselineskip 18pt

\section{Introduction}

In this article, we consider initial-boundary value problems for
a diffusion and a wave equations:
$$\left\{ \begin{array}{rl}
& \ppp_tu(x,t) = \Delta u(x,t) + \mu(t)f(x), \quad x \in \OOO,\,
t>0, \\
& u(x,0) = 0, \quad x \in \OOO,\\
& u(x,t) = 0, \quad x \in \ppp\OOO,\, t>0.
\end{array}\right.
                                \eqno{(1.1)}
$$
$$\left\{ \begin{array}{rl}
& \ppp_t^2u(x,t) = \Delta u(x,t) + \mu(t)f(x), \quad x \in \OOO,\,
t>0, \\
& u(x,0) = \ppp_tu(x,0) = 0, \quad x \in \OOO,\\
& u(x,t) = 0, \quad x \in \ppp\OOO,\, t>0.
\end{array}\right.
                                \eqno{(1.2)}
$$
Here and henceforth $\OOO\subset \R^d$ is a bounded domain with smooth
boundary $\ppp\OOO$ and we set $x=(x_1,..., x_d) \in \R^d$,
$\ppp_j = \frac{\ppp}{\ppp x_j}$, and $\Delta = \sum_{j=1}^d \ppp_j^2$.
Let $\nu = \nu(x)$ be the unit outward normal vector to $\ppp\OOO$ and let
$\ppp_{\nu}u = \nabla u \cdot \nu$.
We mainly consider the zero Dirichlet boundary condition and can treat the
Neumann boundary condition similarly but we omit the details.
Moreover we can consider the inverse problems for (1.1) and (1.2)
where $\Delta$ is replaced
by a suitable elliptic operator with time independent coefficients,
but for simplicity, we mainly argue for $\Delta$.

The source is assumed to be represented in the form of $\mu(t)f(x)$ where
$\mu(t)$ and $f(x)$ describe changes in the time $t$ and the spacial
variable $x$ respectively.  Such a form of separation of variables
is frequently used in modelling diffusion and wave phenomena.
\\

The unique existence of solutions to (1.1) and (1.2) are standard
results (e.g., Evans \cite{E}, Lions and Magenes \cite{LM},
Pazy \cite{Pa}), but we need more regularity of solutions.
We sum up these results as Lemmas 1 and 2.   We arbitrarily fix $T_0>0$.
\\
{\bf Lemma 1.}\\
{\it Let $f\in C^{\infty}_0(\OOO)$ and
$\mu \in H^1(0,T_0)$.
\\
(i) To (1.1), there exists a unique solution $u\in
C([0,T_0]; H^2(\OOO)\cap H^1_0(\OOO)) \cap C^1([0,T_0]; L^2(\OOO))$ and
we can choose a constant $C>0$, dependent on $f$, such that
$$
\Vert u\Vert_{C(\OOO\times [0,T_0])} \le C\Vert \mu\Vert_{L^2(0,T_0)}.
                                        \eqno{(1.3)}
$$
\\
(ii) To (1.2), there exists a unique solution $u\in
C([0, T_0]; H^2(\OOO)) \cap C^1([0, T_0];H^1_0(\OOO))
\cap C^2([0, T_0];L^2(\OOO)) \cap C(\ooo{\OOO} \times [0, T_0])$
such that $\ppp_{\nu}u \in H^1(0,T_0;L^2(\ppp\OOO))$ and (1.3) holds.
}
\\
{\bf Lemma 2.}\\
{\it
Let $f\in L^2(\OOO)$ and $\mu \in C^1[0,T_0]$.
\\
(i) To (1.1), there exists a unique solution $u\in
H^1(0,T_0; L^2(\OOO)) \cap L^2(0,T_0; H^1_0(\OOO))$.
\\
(ii) To (1.2), there exists a unique solution $u\in
H^1(0,T_0; L^2(\OOO)) \cap L^2(0,T_0; H^1_0(\OOO))$ such that
$\ppp_{\nu}u \in H^1(0,T_0;L^2(\ppp\OOO))$.
}

Here we do not aim at the best possible regularity, and
for completeness the proofs of the lemmata are given in Appendix II.
\\

Throughout the article, we assume that $\mu \in L^2_{loc}(\R)$ satisfies
$$
\mu(t) = 0 \quad \mbox{if $t > T$}
$$
with some constant $T>0$.
This means that the a diffusion source for (1.1) and an
external force for (1.2) continue to be activated
only before the moment $T>0$.

In our inverse problems, the measurements can be started after the time
$T>0$, and we are required to determine $\mu(t)$, $0<t<T$ or
$f(x)$, $x\in \OOO$ of the source term.
For example, in the case where the explosion of some equipments such
as nuclear power plant, causes diffusion of contaminants
or dangerous substances, any measurements starting at $t=0$ are not
realistic.   Inverse source problems of determining
$\mu$ or $f$ are well studied if the measurements of data
are started at $t=0$, but to the best knowledge of the authors,
there are no publications on mathematical analysis by data starting
after the time $T>0$.

The main purpose of this article is to establish the uniqueness and
the stability for inverse source problems for (1.1) and (1.2)
by data from the time when the source stopped to be active.

Now we formulate several kinds of inverse source problems and state
our main results.
\\
{\bf \S 1.1. Determination of starting time of decay of source.}

In this subsection, in particular, for $t_0>0$, we consider
$$\left\{ \begin{array}{rl}
& \ppp_tu(x,t) = \Delta u(x,t) + \theta(t-t_0)f(x), \quad x \in \OOO,\,
t>0, \\
& u(x,0) = 0, \quad x \in \OOO,\\
& \ppp_{\nu}u(x,t) = 0, \quad x \in \ppp\OOO,\, t>0,
\end{array}\right.
                                \eqno{(1.4)}
$$
where $\theta\in L^2_{loc}(\R)$ is assumed to be
known and monotone decreasing and satisfy
$$
\theta(t) =
\left\{ \begin{array}{rl}
1, \quad & t \le 0, \\
0, \quad & t\ge a
\end{array}\right.
                                \eqno{(1.5)}
$$
with arbitrarily fixed constant $a>0$.
Only in this subsection, we consider the zero Neumann boundary
condition.
The proof can be modified for the case of the zero
Dirichlet boundary condition.

We note that the diffusion source $\theta(t-t_0)f(x)$
does not act for $t \ge a+t_0$.

Then we consider
\\
{\bf Inverse Problem I.}\\
{\it
Let $f=f(x)$ be known.  Let $T>0$ be sufficiently large and
$x_0\in \OOO$ be arbitrarily chosen.
Determine a starting time $t_0>0$ of the source by
$u(x_0,T)$.}
\\

For known $f$, we assume
$$
f\in C^{\infty}_0(\OOO), \quad f\ge 0, \, \not\equiv 0 \quad
\mbox{on $\ooo{\OOO}$}.                         \eqno{(1.6)}
$$
\\

By $u_{t_0}=u_{t_0}(x,t)$ we note the solution to (1.4), assuming that
$\theta, f$ are fixed.
We are ready to state the uniqueness and the stability as our
main result for Inverse Problem I.\\
{\bf Theorem 1.}\\
{\it
We arbitrarily fix constants $t_*, t^*, a>0$
with $t_* < t^* < T^*-a$.  Then we a priori assume that $t_0, t_1$ are
limited to an interval $(t_*, \, t^*)$:
$$
t_* < t_0, t_1 < t^*.
$$
Then there exists a constant $C=C(t_*, t^*, a,x_0,T^*)>0$ such that
$$
\vert t_1-t_0\vert \le C\vert u_{t_0}(x_0,T^*) - u_{t_1}(x_0,T^*)\vert.
$$
}
\\

This theorem asserts the stability in determining a starting time of decay by
one-shot data $u(x_0, T^*)$, provided that the starting time is assumed to be
in an a priori fixed interval $(t_*, t^*)$.
\\

Now we return to (1.1) and (1.2).
By $u_{\mu}=u_{\mu}(x,t)$ and $u_f=u_f(x,t)$ we denote the solutions to
(1.1) or (1.2) in the cases where we discuss the determination of
$\mu$ and $f$ with fixed $f$ and $\mu$, respectively.
The existence and the regularity of $u_{\mu}$ and $u_f$ are guaranteed by
Lemmata 1 and 2.

Henceforth we will consider the following settings.
First we assume that
$$
\mu(t) = 0\quad \mbox{for $t\ge T$}, \quad T<T_1<T_2.
$$
The exact description of the conditions of $\mu$ and $f$ are different
according to several formulations of our inverse problems, and
later is provided.

We understand that $(T_1, T_2)$ is an observation time interval.
Let $\gamma, \Gamma \subset \ppp\OOO$ be subboundaries.
Now for convenience, we list our main results for the inverse source
problems for (1.1) and (1.2).
\\
{\bf Determination of $\mu(t)$}
\\
{\bf Diffusion equation}
\begin{itemize}
\item
Theorem 2: data $u(x_0,t)$, $T_1<t<T_2$
\item
Proposition 2: data $\ppp_{\nu}u$ on $\gamma \times (T_1,T_2)$
\end{itemize}
{\bf Wave equation}
\begin{itemize}
\item
Theorem 3 (the one-dimensional case), Proposition 1:
data $u(x_0,t)$, $T_1<t<T_2$
\item
Proposition 3: data $\ppp_{\nu}u$ on $\Gamma \times (T_1,T_2)$
\end{itemize}
{\bf Determination of $f(x)$}
\\
{\bf Diffusion equation}
\begin{itemize}
\item
The uniqueness is impossible
in general dimensions by data $u(x_0,t)$, $T_1<t<T_2$.
\item
Theorem 4: data $\ppp_{\nu}u$ on $\gamma \times (T_1,T_2)$
\end{itemize}
{\bf Wave equation}
\begin{itemize}
\item
The uniqueness is impossible in general dimensions
by data $u(x_0,t)$, $T_1<t<T_2$.
\item
Theorem 5: data $\ppp_{\nu}u$ on $\Gamma \times (T_1,T_2)$
\end{itemize}

In general dimensions $d$, the uniqueness does not hold with
data $u(x_0,t)$ for $T_1<t<T_2$, because
unknown $f$ depends on $d$-variables, but the data
depend only on one variable $t$.
\\

Here we do not discuss inverse problems with final data $u(\cdot,T)$,
and for the heat equations, we refer to Cheng and Liu \cite{CJ},
Choulli and Yamamoto \cite{ChY} and the references therein.

{\bf \S1.2. Determination of $\mu(t)$ of the source term by pointwise data}

In this subsection, we assume:
$$
\mu \in H^1_{loc}(0, \infty), \quad \not\equiv 0, \quad
\mu(t) = 0\quad \mbox{for $t\ge T$}                         \eqno{(1.7)}
$$
and
$$
f \in C^{\infty}_0(\OOO), \quad f \not\equiv 0.      \eqno{(1.8)}
$$
Here we set
$$
H^1_{loc}(0,\infty):= \{ \mu;\,
\mu\vert_{(0,T_0)} \in H^1(0,T_0) \quad \mbox{with any $T_0>0$}\}.
$$
We can relax the regularity of $f$ but we assume (1.8) for simplicity.

Moreover, let $T_1, T_2$ be given such that
$$
T<T_1<T_2
$$
and let $x_0\in \OOO$ be arbitrarily chosen.

We consider

{\bf Inverse Problem II}.\\
{\it
In (1.1) and (1.2), we are given $f$. Determine $\mu(t)$, $0<t<T$, by
$u(x_0,t)$ for $T_1<t<T_2$.
}
\\

By Lemma 1, we know that $u\in C(\ooo{\OOO} \times [T_1, T_2])$, and so
our observation data $u(x_0,t)$, $T_1<t<T_2$, are well-defined.

For the statement of our main results for Inverse Problem II,
we introduce notations.
Let $\la_j$, $j\in \N$ be all the distinct eigenvalues of the operator
$A=-\Delta$ with the domain $\mathcal{D}(A) = H^2(\OOO) \cap H^1_0(\OOO)$.
We know that $\la_j > 0$ for all $j\in \N$.
By $d_j$ we denote the multiplicity of $\la_j$, $j\in \N$, and
$\{ \va_{jk}\}_{1\le k\le d_j}$ is an orthonormal basis of
Ker $(\la_j-A) := \{ v\in \mathcal{D}(A);\, Av = \la_jv\}$,
and we define the orthogonal projection $P_j$ from
$L^2(\OOO)$ to \\
Ker $(\la_j-A)$ by
$$
P_jf = \sum_{k=1}^{d_j} (f, \va_{jk})\va_{jk}.           \eqno{(1.9)}
$$
We set
$$
\Lambda = \Lambda(x_0) := \{ j\in \N;\, (P_jf)(x_0) = 0\}.
                                                         \eqno{(1.10)}
$$

We are ready to state our main result for Inverse Problem II for
the diffusion equation (1.1).
\\
{\bf Theorem 2.}\\
{\it
Assume that (1.8) holds and $\mu_1, \mu_2$ satisfy (1.7).
Let $x_0 \in \OOO$ be arbitrarily chosen.  Then
$$
u_{\mu_1}(x_0,t) = u_{\mu_2}(x_0,t) \quad \mbox{for $T_1<t<T_2$}
                            \eqno{(1.11)}
$$
yields $\mu_1(t) = \mu_2(t)$ for $0<t<T$ if and only if
$$
\sum_{j\in \N \setminus \Lambda(x_0)} \frac{1}{\la_j} = \infty.
                                                \eqno{(1.12)}
$$
}

{\bf Corollary.}\\
{\it
The uniqueness holds for Inverse Problem II only if
the spatial dimensions $d\ge 2$.
The uniqueness always fails for $d=1$.
}
\\

In the case of $d=1$, our data $u_{\mu}(x_0,t)$, $t>T$ starting after the
incident (i.e., $\mu(t) = 0$ for $t>T$), cannot give the uniqueness.
Moreover, as is seen by the proof in Section 3, even if we will take the
perfect observation data $u_{\mu}(x,t)$
with all $x\in \OOO$ and $t>T$, we can determine at
best
$$
\int^T_0 e^{\la_ns}\mu(s) ds, \quad n\in \N,
$$
which is the same information determined by currently adopted
pointwise data $u_{\mu}(x_0,t)$, $t<T$.
In other words, even the perfect observation data cannot augment
any information of poinwise data $u(x_0,t)$, $t>T$ with fixed $x_0
\in \OOO$.
\\
\vspace{0.2cm}
{\bf Example of (1.12).}\\
Let $d=2$ and $\OOO = (0, \ell_1) \times (0, \ell_2)$ with $\ell_1, \ell_2
> 0$.  Then we can directly verify that
$$
\{ \la_j\}_{j\in \N} = \left\{
\left( \frac{m_1^2}{\ell_1^2} + \frac{m_2^2}{\ell_2^2}\right)\pi^2;\,
m_1, m_2 \in \N\right\}.            \eqno{(1.13)}
$$
For simplicity, we consider only the case where $\frac{\ell_1^2}{\ell_2^2}$
is an irrational number.
Then $\frac{m_1^2}{\ell_1^2} + \frac{m_2^2}{\ell_2^2}
= \frac{n_1^2}{\ell_1^2} + \frac{n_2^2}{\ell_2^2}$ with $m_1, m_2, n_1, n_2
\in \N$ imply $m_1=n_1$ and $m_2=n_2$, which means that the multiplicity
of $\la_j$, $j\in \N$ is one.  We re-number
$\left(\frac{m_1^2}{\ell_1^2} + \frac{m_2^2}{\ell_2^2}\right)\pi^2$
with $m_1, m_2 \in \N$ as $0<\la_1 < \la_2 < \cdots$.
Then for $\la_j = \left( \frac{m_1^2}{\ell_1^2} + \frac{m_2^2}{\ell_2^2}
\right)\pi^2$, we choose
$$
\va_j(x) = \frac{2}{\sqrt{\ell_1\ell_2}}\sin \frac{m_1\pi x_1}{\ell_1}
\sin \frac{m_2\pi x_2}{\ell_2}.
$$
We choose a monitoring point $x_0 = (x_0^1, x_0^2)
\in (0, \ell_1) \times (0,\ell_2)$ such that
$$
\mbox{$\frac{x_0^1}{\ell_1}$ and $\frac{x_0^2}{\ell_2}$ are
irrational numbers.}                        \eqno{(1.14)}
$$
Then $m_1\frac{x_0^1}{\ell_1}, m_2\frac{x_0^2}{\ell_2}
\not\in \N$ for all $m_1, m_2 \in \N$, and so $\va_j(x_0^1,x_0^2) \ne 0$.
Thus, $(P_jf)(x_0) = 0$ if and only if
$\int_{\OOO} \va_j(x)f(x) dx = 0$, and so we see that
$\int_{\OOO} \va_k(x)f(x) dx \ne 0$
except for a finite number of $k\in \N$, then
we can conclude that $\Lambda(x_0)$ is a finite set.

Under (1.14), we can verify that (1.12) is satisfied.
Indeed, by Agmon \cite{Ag} or Courant and Hibert \cite{CH}
for example, we know that $\la_j = \rho_0j + o(1)$ as $j \to \infty$, where
$\rho_0>0$ is a constant.  Therefore, $\sum_{j=1}^\infty \frac{1}{\la_j}
= \infty$.  Thus if $\frac{\ell_1^2}{\ell_2^2}$ is an irrational number and
$\int_{\OOO} \va_k(x)f(x) dx \ne 0$ except for a finite number of
$k\in \N$, then (1.12) is satisfied.
\\

In the case where the measurement starts at $t=0$, that is, $T_1=0$, the
inverse source problem is even stable (e.g., Cannon and Esteva \cite{CE},
Saitoh, Tuan and Yamamoto \cite{STY}).
However, for $T_1>0$, to the best knowledge of the authors, there have been
no publications, although such formulated inverse problems are practical.

We can expect to establish conditional stability which holds under
suitable a priori boundedness condition on $\mu(t)$, and we
conjecture that the rate is of weak type such as logarithmic rate.

Moreover, in the case where $0<T_1<T_2<T$, that is, the observation
starts after the activities of the source but before the stop of the
activities, the uniqueness seems impossible, but here we do not discuss
the details.
\\

Next we consider the wave equation (1.2).  For the uniqueness,
we can prove positive results only for the one-dimensional case
$d =1$.
\\
{\bf Theorem 3 (uniqueness for time dependent factor of wave source).}\\
{\it
We assume (1.8) and (1.7) for $\mu_1, \mu_2$.
Let $d=1$ and $\OOO = (0, \ell)$.  Then $u_{\mu_1}(x_0,t)
= u_{\mu_2}(x_0,t)$ for $T_1<t<T_2$ yields $\mu_1(t) = \mu_2(t)$ for
$0\le t\le T$ if and only
$$
\left\{ \begin{array}{rl}
& \frac{x_0}{\ell} \in \R \setminus \Q, \quad T\le 2\ell, \quad
T_2-T_1 \ge 2\ell, \\
& \int^{\ell}_0 f(x) \sin \frac{n\pi}{\ell}x dx \ne 0,
\quad n\in \N.
\end{array}\right.
                                    \eqno{(1.15)}
$$
}
\\

The assumption $T_2-T_1 \ge 2\ell$ in (1.15) means that
we have to take measurements longer than
$2\ell$, while the width of support of $\mu(t)$ should not be long, that is,
$\le 2\ell$.
We do not know the uniqueness for general dimensions $d\ge 2$.

For the case of $d\ge 2$, we can prove the uniqueness with different
measurements of pointwise data $u(x_0,t)$.
We recall that $\Lambda(x_0) \subset \N$ is defined by (1.10).

We assume that $\N \setminus \Lambda(x_0)$ is an infinite subset of $\N$
and
$$
\lim_{k\to \infty, k\in \N \setminus \Lambda(x_0)}
\frac{k}{\sqrt{\la_k}} > 0.
                                                 \eqno{(1.16)}
$$
We remark that also the existence of the limit is assumed in (1.16).
\\
{\bf Proposition 1.}\\
{\it
Let $d\ge 2$ and $f \in C^{\infty}_0(\OOO)$, and
$$
T < \pi\lim_{k\to \infty, k\in \N \setminus \Lambda(x_0)}
\frac{k}{\sqrt{\la_k}}.                                         \eqno{(1.17)}
$$
Then
$$
\limsup_{t\to \infty} \vert u_{\mu_1}(x_0,t) - u_{\mu_2}(x_0,t)\vert = 0
                                                   \eqno{(1.18)}
$$
yields $\mu_1(t) = \mu_2(t)$ for $0\le t\le T$.
}
\\

If $\lim_{k\to \infty, k\in \N \setminus \Lambda(x_0)}
\frac{k}{\sqrt{\la_k}} = \infty$, then we interpret that (1.17) holds true for
any $T> 0$.

The observation in (1.18) is concerned with the asymptotics of
$u_{\mu}(x_0,t)$ as $t \to \infty$ and (1.18) requires that the
data at $x_0$ are asymptotically equal as $t \to \infty$.

We conclude this subsection with
\\
{\bf Remark on condition (1.16).}
\\
We consider one sufficien condition for (1.16) in terms of the multiplicities
of $\la_k$.
We assume
$\Lambda(x_0) := \{ j\in \N;\, (P_jf)(x_0) = 0\}
= \emptyset$.
We recall that by $\la_1 < \la_2 < \cdots$ we number the set of all the
eigenvalues, not taking the multiplicities into consideration.
On the other hand, by $\sigma_k$, $k\in \N$,
we number all the eigenvalues according to their multiplicities:
$0<\sigma_1 \le \sigma_2 \le \sigma_3 \le \cdots \longrightarrow \infty$.
In other words, the value $\lambda_j$ appears $d_j$-times among
the sequence $\sigma_1, \sigma_2, \sigma_3, \cdots$, where
$d_j:= \mbox{dim Ker}\,(\la_j-A)$.
More precisely, setting
$$
m_k := \sum_{i=1}^k d_i, \quad k\in \N,
$$
we see
$$
\sigma_j =
\left\{ \begin{array}{rl}
\la_1, \quad &1\le j \le m_1, \\
\la_2, \quad &m_1+1\le j \le m_2, \\
\cdots & \cdots,\\
\la_k, \quad & m_{k-1}+1\le j \le m_k.
\end{array}\right.
$$

By the definition, we note
$$
\sigma_j \le \la_j, \quad j\in \N.
$$
As general information of
$\sigma_j$ the following is classical:
$$
\sigma_j = \rho_0 j^{\frac{2}{d}} + o(d^{\frac{2}{d}}) \quad
\mbox{as $j \to \infty$}                       \eqno{(1.19)}
$$
(e.g., Theorems 14.6 and 15.1 in Agmon \cite{Ag}, Chapter 6 in
Courant and Hilbert \cite{CH}).  Here $\rho_0>0$ is a constant
determined by $d$ and $\vert \OOO\vert$.

By (1.19) we have
$$
\rho_0 m_{k-1}^{\frac{2}{d}}(1+o(1))
\le \la_k \le \rho_0 m_k^{\frac{2}{d}}(1+o(1))
$$
as $k\to \infty$.  We here note that if $j\to \infty$, then
$k\to \infty$.
Therefore
$$
\lim_{k\to \infty} \frac{k}{m_k^{\frac{1}{d}}} > 0
                                                         \eqno{(1.20)}
$$
yields (1.16).  Indeed
$$
\frac{k}{\sqrt{\la_k}} \ge \frac{k}{m_k^{\frac{1}{d}}}
\frac{1}{\sqrt{\rho_0}(1+o(1))}
$$
as $k \to \infty$.

In the following two examples, we consider (1.20).
\\
{\bf Example 1: Let $d=2$ and $\OOO = (0, \ell_1) \times (0, \ell_2)$.}
\\
First we assume $\frac{\ell_1^2}{\ell_2^2} \not\in \Q$.  By (1.13),
all the eigenvalues are simple, that is, $d_i=1$ for $i\in \N$.
Therefore $m_k = k$ for $k\in \N$, that is,
$$
\lim_{k\to \infty} \frac{k}{\sqrt{m_k}} = \infty.
$$
Consequently, we need not any assumption for $T$.
Also in general domain $\OOO$, if all the eigenvalues are simple
except for a finite number, then (1.16) is true in terms of (1.19), and
(1.17) is satisfied for any $T>0$.

Next we assume $\frac{\ell_1^2}{\ell_2^2} \in \Q$.
Then all the eigenvalues are not necessarily simple, and we do not know
suitable estimates of $m(k)$ as $k \to \infty$, but we have
$$
\limsup_{i\to \infty} d_i = \infty,
$$
whose proof is found in Yamamoto \cite{Ya4} for example.
We do not know how rapidly $d_i$ tends to $\infty$.  If they goes to
$\infty$ very rapidly, then
$\inf\lim_{k\to\infty} \frac{k}{\sqrt{m_k}} = 0$ may happen, which
breaks (1.20), and we do not know whether (1.16) holds.
\\
{\bf Example 2: Let $d=2$ and $\OOO = \{ x\in \R^2;\, \vert x\vert < 1\}$.}
\\
It is known that
$d_i = \left\{ \begin{array}{rl}
1, \quad \mbox{if $i=1$}, \\
2, \quad \mbox{if $i\ge 2$}.
\end{array}\right.$.
Therefore $m_k = 2k-1$ for $k\in \N$ and (1.16) holds.
We see that (1.17) holds for any $T>0$.
\\
\vspace{0.2cm}
\\
{\bf \S 1.3. Determination of $\mu(t)$ of the source term by
boundary data.}

In this subsection, for the inverse problems, we adopt data $\ppp_{\nu}u$
on a lateral subboundary.

We set
$$
\www{\Lambda} := \{ j\in \N;\, P_jf = 0 \quad
\mbox{in $\OOO$} \}.                       \eqno{(1.21)}
$$
For the diffusion equation, we can show
\\
{\bf Proposition 2.}\\
{\it
Let $\gamma \subset \ppp\OOO$ be an arbitrarily chosen subboundary.
We assume (1.7) and (1.8).  Then $\ppp_{\nu}u_{\mu_1} =
\ppp_{\nu}u_{\mu_2}$ on $\gamma \times (T_1,T_2)$ yields
$\mu_1(t) = \mu_2(t)$ for $0<t<T$ if and only if
$$
\sum_{j\in \N\setminus \www{\Lambda}} \frac{1}{\la_j} = \infty.
                                                        \eqno{(1.22)}
$$
}

The inverse source problem with boundary data is overdetermining because
an unknown function depends only on one variable $t$.
The condition (1.22) for the uniqueness is weaker than (1.12) for the
uniqueness by the pointwise data $u(x_0,t)$, $T_1<t<T_2$.
\\

In the one-dimensional case, Theorem 3 asserts the uniqueness for a
wave equation by the pointwise data $u(x_0,t)$ for $T_1<t<T_2$,
but we do not know the corresponding uniqueness for general dimensions.
On the other hand, by boundary observations, we can conclude
the uniqueness for general
dimensions $d \ge 2$, as the following proposition shows.
\\
{\bf Proposition 3.}\\
{\it
Let
$$
\lim_{k\to \infty,k\in \N\setminus \www{\Lambda}}
\frac{k}{\sqrt{\la_k}} > 0.    \eqno{(1.23)}
$$
We assume that there exists $x_0 \in \R^d$ such that
$$
\Gamma \supset \{ x\in \ppp\OOO;\, ((x-x_0)\cdot \nu(x)) \ge 0\}, \quad
T_2-T_1 > 2\max_{x\in \ooo{\OOO}} \vert x-x_0\vert
                                                            \eqno{(1.24)}
$$
and
$$
T < \pi \lim_{k\to \infty, k\in \N\setminus \www{\Lambda}}
\frac{k}{\sqrt{\la_k}}.             \eqno{(1.25)}
$$
Moreover, let $\mu_1, \mu_2 \in C^1[0,T]$ satisfy
$\mu_1(t) = \mu_2(t) = 0$ for $t > T$.
Then $\ppp_{\nu}u_{\mu_1}(x,t) = \ppp_{\nu}u_{\mu_2}(x,t) = 0$ for
$x \in \Gamma$ and $T_1<t<T_2$ implies $\mu_1(t) = \mu_2(t)$ for $0<t<T$.
}
\\

The condition (1.25) can be characterized similarly to (1.20).
The observation time length $T_2-T_1$ should be
longer, while the duration time $T$ of the source needs
to be bounded for the uniqueness for our inverse source problem.
\\
\vspace{0.2cm}
{\bf \S 1.4. Determination of $f(x)$ of source terms by boundary data.}

In this subsection, we assume that $\Gamma, \gamma \subset \ppp\OOO$ are
subboundaries, and $T<T_1<T_2$,
$$
\mu \in C^1[0,\infty), \quad \mu(t) = 0 \quad
\mbox{if $t\ge T$}                               \eqno{(1.26)}
$$
and
$$
f \in L^2(\OOO).                                 \eqno{(1.27)}
$$

We consider
\\
{\bf Inverse Problem III.}\\
{\it
Determine $f=f(x)$, $x\in \OOO$ by
$\ppp_{\nu}u$ on some subboundary over a time interval
$(T_1, T_2)$.
}
\\

In the case of $T_1=0$, there are many researches on the uniqueness and the
stability.  Here we limit ourselves to a few works: for the diffusion equation
(1.1), we refer to Xu and Yamamoto \cite{XY}, Yamamoto \cite{Ya1}, while
we can consult Yamamoto \cite{Ya2}, \cite{Ya3} for the wave equation (1.2).

First we consider the determination of $f(x)$ for (1.1).
Henceforth for arbitrarily chosen constants $M>0$ and $\ell \in \N$,
we define an admissible set of unknowns $f$ by
$$
\mathcal{F}_{M,\ell} := \{ f\in \mathcal{D}(A^{\ell});\,
\Vert A^{\ell}f\Vert_{L^2(\OOO)} \le M\}.                \eqno{(1.28)}
$$
\\
{\bf Theorem 4 (determination of $f(x)$ in the diffusion equation).}\\
{\it
Let $\gamma$ be an arbitrarily chosen subboundary of $\ppp\OOO$.
\\
(i) We assume
$$
\int^T_0 e^{\la_kt}\mu(t) dt \ne 0, \quad k \in \N.
                                                 \eqno{(1.29)}
$$
Then $\ppp_{\nu}u = 0$ on $\gamma \times (T_1,T_2)$ yields
$f=0$ in $\OOO$.
\\
(ii) We further assume
$$
\mu \ge 0, \, \not\equiv 0 \quad \mbox{on $[0,T]$}.      \eqno{(1.30)}
$$
Then for $0< \theta < \ell$, we can find a constant $C=C(M,\theta)>0$
such that
$$
\Vert f\Vert_{L^2(\OOO)} \le C\left(
\frac{1}{\log \frac{1}{\Vert \ppp_{\nu}u\Vert_{L^2(\gamma \times (T_1,T_2))}}}
\right)^{\theta}
$$
as $\Vert \ppp_{\nu}u\Vert_{L^2(\gamma \times (T_1,T_2))} \longrightarrow 0$
for each $f \in \mathcal{F}_{M,\ell}$.
}
\\

In the case of $T_1=0$, we have the same estimate without assumption
(1.29) (e.g., Yamamoto \cite{Ya1}).
\\
{\bf Remark.}
We need the assumption (1.29).  Indeed, let $\va$ be an eigenfunction for
$\la_1$ of $-\Delta$ with the zero Dirichlet boundary condition and
let (1.29) fail, for example, let $\int^T_0 e^{\la_1s} \mu(s) ds = 0$.
Then, we readily see that
$$
u(x,t):= \left( \int^t_0 e^{-\la_1(t-s)}\mu(s) ds\right) \va_1(x),
\quad x \in \OOO, \, t>0,
$$
satisfies
$$
\left\{ \begin{array}{rl}
& \ppp_tu = \Delta u(x,t) + \mu(t)\va_1(x), \quad x\in \OOO, \, t>0,\\
& u(x,t) = 0, \quad x\in \ppp\OOO, \, t>0, \\
& u(x,0) = 0, \quad x \in \OOO,
\end{array}\right.
$$
and
$$
u(x,t) = e^{-\la_1t} \left( \int^T_0 e^{\la_1s}\mu(s) ds\right) \va_1(x)
= 0, \quad x \in \OOO, \, t>T,
$$
because $\mu(s) = 0$ for $s\ge T$ and $\int^T_0 e^{\la_1s}\mu(s) ds = 0$.

Therefore even $u = 0$ on $\ooo{\OOO} \times (T_1,T_2)$ holds, and
in particular, $\ppp_{\nu}u\vert_{\ppp\OOO\times (T_1,T_2)} = 0$.
However, $f=\va_1 \ne 0$, that is, the uniqueness fails without (1.29).
\\

Next we show
\\
{\bf Theorem 5 (determination of $f(x)$ in the wave equation).}\\
{\it
We assume (1.24) and
$$
\mu(0) \ne 0
$$
and
$$
\left\vert \int^T_0 \mu(s) \sin \sqrt{\la_n}(T-s) ds \right\vert
+ \left\vert \int^T_0 \mu(s) \cos \sqrt{\la_n}(T-s) ds \right\vert
\ne 0 \quad \mbox{for all $n\in \N$}.       \eqno{(1.31)}
$$
Then there exists a constant $C>0$ such that
$$
C^{-1}\Vert \ppp_{\nu}u\Vert_{H^1(T_1,T_2;L^2(\ppp\OOO))}
\le \Vert f\Vert_{L^2(\OOO)}
\le C\Vert \ppp_{\nu}u\Vert_{H^1(T_1,T_2;L^2(\Gamma))}
                                                           \eqno{(1.32)}
$$
for each $f \in L^2(\OOO)$.
}
\\

In the case of $T_1=0$, we refer to the existing works: for example,
Yamamoto \cite{Ya2}, \cite{Ya3} which proves (1.32) only by (1.24) and
$\mu(0) \ne 0$, not assuming (1.31).
\\
{\bf Example of (1.31).}\\
We verify that $\mu(s) = T-s$ satisfies (1.31).
Indeed we directly have
$$
\int^T_0 (T-s)\sin (T-s)\sqrt{\la_n} ds
= \frac{-T\sqrt{\la_n}\cos T\sqrt{\la_n} + \sin T\sqrt{\la_n}}{\la_n}
$$
and
$$
\int^T_0 (T-s)\cos (T-s)\sqrt{\la_n} ds
= \frac{T\sqrt{\la_n}\sin T\sqrt{\la_n} + \cos T\sqrt{\la_n} - 1}{\la_n}.
$$
Assume that both are zero for some $n_0\in \N$.
Setting $\xi = T\sqrt{\la_{n_0}} \ne 0$, we see that
$$
\xi \cos \xi - \sin \xi = 0, \quad \cos\xi + \xi\sin \xi = 1.
$$
Hence $\cos\xi = \frac{1}{1+\xi^2}$ and $\sin \xi = \frac{\xi}{1+\xi^2}$,
By $\cos^2\xi + \sin^2\xi = 1$, we obtain $\frac{1}{1+\xi^2}=1$, which yields
$\xi=0$, which is a contradiction by $\xi = \sqrt{\la_{n_0}}T \ne 0$.
Thus $\mu(s) = T-s$ satisfies (1.31).
\\

This article is composed of 9 sections.  In Sections 2 and 3, we prove
Theorems 1 and 2 respectively.  In Section 4, we prove Theorem 3 and
Proposition 1.  Sections 5 and 6 are devoted to the proofs of
Theorems 4 and 5.  In Section 7, we prove Propositions 2 and 3.
For convenience, Section 8 provides proofs of standard
uniqueness results for the case of $T_1=0$ in determining $\mu(t)$.
Section 9 gives the proofs of Lemmata 1 and 2.
\section{Proof of Theorem 1}

Without loss of generality, we can assume
that $0<t_0<t_1$.  For $t\ge 0$, we set
$$
g(t) := \theta(t-t_1) - \theta(t-t_0) = 0
\quad \mbox{if $t \le t_0$ or $t>t_1 + a$.}
                                   \eqno{(2.1)}
$$

Setting $y=u_{t_1} - u_{t_0}$, we have
$$
\left\{ \begin{array}{rl}
& \ppp_ty(x,t) = \Delta y(x,t) + g(t)f(x), \quad x\in \OOO, \, t>0,\\
& \ppp_{\nu}y(x,t) = 0, \quad x\in \ppp\OOO, \, t>0, \\
& y(x,0) = 0, \quad x \in \OOO.
\end{array}\right.                       \eqno{(2.2)}
$$
Let $G = G(x,y,s)$ be the Green function for $\ppp_t - \Delta$
in $\OOO$ with the
zero Neumann boundary condition.  Then
$$
G(x,y,s) > 0, \quad x,y \in \OOO, \, 0<s<T
$$
(e.g., Theorem 10.1 (i) in It{\^o} \cite{I}).
Setting $D:= \mbox{supp}\, f$, by (1.6) we see that
$D \subset \OOO$.  Note that $D$ is a closed set and
$T^* > t^* + a$ by the assumption.
Since $G(x_0,y,s)$ is a continuous function in $y \in D$ and
$T^* - t^* -a \le s \le T^*-t_*$,
we can choose a constant $C_0 = C_0(D,t_*,t^*+a)>0$
such that
$$
G(x_0,y,s) \ge C_0 > 0, \quad y\in D, \,
T^* - t^* -a \le s \le T^*-t.
                                   \eqno{(2.3)}
$$
Moreover, we have
$$
y(x,t) = \int^t_0 \int_{\OOO} G(x,y,t-s) f(y)g(s) dy ds,
\quad x \in \OOO, \, t>0                          \eqno{(2.4)}
$$
(e.g., \cite{I}).
In general, (2.4) holds for $g \in C[0,T^*]$, and we here do not assume that
$g$ is continuous.  However, (1.3) in Lemma 1 asserts
$$
\Vert y\Vert_{C(\ooo{\OOO}\times [0,T^*])} \le C\Vert g\Vert_{L^2(0,T^*)}
                                                    \eqno{(2.5)}
$$
in view of $f \in C^{\infty}_0(\OOO)$.  Therefore by taking an
approximating sequence $g_n \in C^{\infty}_0(0,T^*)$ to $g$ in $L^2(0,T^*)$,
we can verify that (2.4) holds for $g \in L^2(0,T^*)$.

By $0<t_* < t^*+a < T^*$, (1.6), and (2.3), we obtain
\begin{align*}
& y(x_0,T^*) = \int^{T^*}_0 \int_{\OOO}
G(x_0,y,T^*-s)f(y) g(s) dy ds\\
\ge & \int^{t^*+a}_{t_*} \int_D G(x_0,y,T^*-s)f(y) g(s) dy ds\\
\ge& C_0\left(\int^{t^*+a}_{t_*} g(s) ds\right)
\int_D f(y) dy
\ge C_0C_1 \int^{t^*+a}_{t_*} g(s) ds .
\end{align*}
Here, by (1.6), we used $C_1 := \int_D f(y) dy > 0$.
On the other hand,
\begin{align*}
& \int^{t^*+a}_{t_*} g(s) ds
= \int^{t^*+a}_{t_*} \theta(s-t_1) ds
- \int^{t^*+a}_{t_*} \theta(s-t_0) ds\\
=& \int^{a+t^*-t_1}_{t_*-t_1} \theta(\eta) d\eta
- \int^{a+t^*-t_0}_{t_*-t_0} \theta(\eta) d\eta\\
=& \left( \int^0_{t_*-t_1} + \int^a_0 + \int^{a+t^*-t_1}_a\right)
\theta(\eta) d\eta
- \left( \int^0_{t_*-t_0} + \int^a_0 + \int^{a+t^*-t_0}_a\right)
\theta(\eta) d\eta
= \int^{t_*-t_0}_{t_*-t_1} \theta(\eta) d\eta.
\end{align*}
Here we used $\theta(\eta) = 0$ for $\eta \ge a$, and noted that
$a<a+t^*-t_1$ and $a<a+t^*-t_0$.  By $t_*-t_0, t_*-t_1 < 0$, we obtain
$\int^{t^*+a}_{t_*} g(s) ds = \int^{t_*-t_0}_{t_*-t_1} d\eta = t_1-t_0$.

Noting that $y = u_{t_1} - u_{t_0}$, we complete the proof of Theorem 1.
\section{Proof of Theorem 2}

Let $\mu = \mu_1 - \mu_2$ and $u=u_{\mu_1} - u_{\mu_2}$.  By the
eigenfunction expansion, we have
$$
u(x,t) = \sumn \left(\int^t_0 e^{-\la_n(t-s)}\mu(s) ds\right) P_nf(x), \quad
x\in \OOO, \, t>0                         \eqno{(3.1)}
$$
(e.g., \cite{I}).  In terms of $f \in C^{\infty}_0(\OOO)$ and
Lemma 1 in Section 1, we can verify that the series in (3.1) is
convergent in $C(\ooo{\OOO}\times [0,T_2])$.

Let $t \in [T_1, T_2]$.  Since $\mu(s) = 0$ for $s>T$, we have
$$
u(x_0,t) = \sumn e^{-\la_nt}\left( \int^T_0 e^{\la_ns}\mu(s) ds\right)
P_nf(x_0)
$$
$$
= \sum_{n\in \N \setminus \Lambda(x_0)} e^{-\la_nt}\left( \int^T_0
e^{\la_ns}\mu(s) ds\right)(P_nf)(x_0) = 0, \quad t>T.             \eqno{(3.2)}
$$
Here we recall that the set $\Lambda(x_0) \subset \N$ is defined by (1.10).

We know
\\
{\bf Lemma 3.}\\
{\it
Let $\sumn \vert a_n\vert e^{-\la_nT_1} < \infty$.  If
$$
\sumn a_n e^{-\la_nt} = 0, \quad T_1<t<T_2,
$$
then $a_n=0$ for all $n\in \N$.
}
\\
{\bf Proof of Lemma 3.}\\
For completeness, we provide the proof although the lemma is well-known.
We set $b_n:= a_ne^{-\la_nT_1}$.
By
$$
\sumn \vert a_n\vert e^{-\la_nT_1} = \sumn \vert b_n\vert < \infty,
$$
we see that
the series $\sumn b_ne^{-\la_n(z-T_1)}$ converges uniformly and absolutely
in any
compact domain in $\{z\in \C;\, \mbox{Re}\, z > T_1\}$. Therefore
$h(z) := \sumn b_n e^{-\la_n(z-T_1)}$ is analytic if Re $z > T_1$.
Hence $\sumn b_ne^{-\la_n(t-T_1)} = 0$ for all $t>T_1$ by the analyticity.
Since $\la_1 < \la_2 < \cdots$, we have
$$
b_1 + \sum_{n=2}^{\infty} b_ne^{-(\la_n-\la_1)(t-T_1)} = 0
\quad \mbox{for $t>T_1$}.
$$
Since
$$
\left\vert \sum_{n=2}^{\infty} b_ne^{-(\la_n-\la_1)(t-T_1)} \right\vert
\le e^{-(\la_2-\la_1)(t-T_1)} \sum_{n=2}^{\infty} \vert b_n \vert
$$
and
$$
\sum_{n=2}^{\infty} \vert b_n\vert
= \sum_{n=2}^{\infty} \vert a_n\vert e^{-\la_nT_1} < \infty,
$$
letting $t \to \infty$, we see that $b_1=0$, that is, $a_1=0$.
Continuing this argument,
we reach $a_n=0$ for $n\in \N$.  Thus the proof of Lemma 3 is
complete.  $\blacksquare$
\\

We return to the proof of Theorem 2.
Setting
$$
a_n = \left(\int^T_0 e^{\la_ns} \mu(s) ds\right) P_nf(x_0), \quad n \in \N,
$$
by (3.2) we obtain
$$
\sum_{n\in \N\setminus \Lambda(x_0)} a_ne^{-\la_nt} = 0, \quad T_1<t<T_2
$$
and
\begin{align*}
& \sum_{n\in \N\setminus \Lambda(x_0)} \vert a_ne^{-\la_nt}\vert
= \sumn \left\vert \left(\int^T_0 e^{-\la_n(t-s)}\mu(s) ds \right)
(P_nf)(x_0) \right\vert\\
\le& \sumn \left(\int^T_0 \vert \mu(s)\vert ds\right) \vert (P_nf)(x_0) \vert,
\quad t\ge T_1 > T.
\end{align*}
We recall that $A = -\Delta$ with $\mathcal{D}(A) = H^2(\OOO) \cap
H^1_0(\OOO)$.  By (1.8) we see that $f \in \mathcal{D}(A^m)$ with
each $m\in \N$.
Since $P_jf \in \mathcal{D}(A^m)$ and $AP_jf = \la_jP_jf$, choosing
even number $\beta > \frac{d}{2}$, by Sobolev embedding, we obtain
$$
\vert P_jf(x_0)\vert \le \Vert P_jf\Vert_{C(\ooo{\OOO})}
\le C\Vert P_jf\Vert_{H^{\beta}(\OOO)}
\le C\Vert A^{\frac{\beta}{2}}P_jf\Vert_{L^2(\OOO)}
= C\la_j^{\frac{\beta}{2}}\Vert P_jf\Vert_{L^2(\OOO)}
$$
and
$$
\Vert P_jf\Vert^2_{L^2(\OOO)} = \sum_{k=1}^{d_j} \vert (P_jf, \va_{jk})\vert^2
= \sum_{k=1}^{d_j} \vert (A^mP_jf, \, A^{-m}\va_{jk})\vert^2
= \la_j^{-2m}\sum_{k=1}^{d_j} \vert (A^mP_jf, \va_{jk})\vert^2.
$$
Therefore,
\begin{align*}
& \vert P_jf(x_0)\vert^2
\le C\la_j^{-2m+\beta}\sum_{k=1}^{d_j} \vert (A^mP_jf, \va_{jk})\vert^2\\
\le& C\Vert A^mP_jf\Vert^2_{L^2(\OOO)}\la_j^{-2m+\beta}
\le C\Vert A^mf\Vert^2_{L^2(\OOO)}\la_j^{-2m+\beta}, \quad j\in \N.
\end{align*}
Moreover, we choose
$$
\beta_1 > \frac{d}{4} \quad\mbox{and}\quad
m > \frac{1}{2}\left( \beta + 2\beta_1 + \frac{d}{2} \right),
$$
and we re-number the eigenvalues
$\sigma_j$ according to the multiplicities: $0 < \sigma_1 \le \sigma_2 \le
\sigma_3 \le \cdots \longrightarrow$ and we note that
$\la_j\ge \sigma_j$ by the re-numbering.
Then $\la_j \ge \sigma_j \sim \rho_0j^{\frac{2}{d}}$ as $j \to \infty$
(e.g., \cite{Ag}, \cite{CH}) and
\begin{align*}
& \sum^{\infty}_{j=1} \vert P_jf(x_0)\vert
= \sum_{j=1}^{\infty} \la_j^{-\beta_1}\la_j^{\beta_1}\vert P_jf(x_0)\vert
\le \left( \sumj \la_j^{-2\beta_1}\right)^{\frac{1}{2}}
\left( \sumj \la_j^{2\beta_1} \vert P_jf(x_0)\vert^2 \right)^{\frac{1}{2}}\\
\le& C\Vert A^mf\Vert_{L^2(\OOO)}
\left( \sumj \la_j^{-2\beta_1}\right)^{\frac{1}{2}}
\left( \sumj \la_j^{-2m+\beta+2\beta_1} \right)^{\frac{1}{2}}\\
\le &C\Vert A^mf\Vert_{L^2(\OOO)}
\left( \sumj j^{\frac{-4\beta_1}{d}} \right)^{\frac{1}{2}}
\left( \sumj j^{\frac{2}{d}(-2m+\beta+2\beta_1)} \right)^{\frac{1}{2}}.
\end{align*}
Here $m > \frac{1}{2}\left(\beta+2\beta_1+\frac{d}{2}\right)$
and $\beta_1 > \frac{d}{4}$ yield that
$\frac{2}{d}(-2m+\beta+2\beta_1) < -1$ and
$-\frac{4\beta_1}{d}< -1$, so that
we see that the right-hand side is convergent, that is,
$$
\sumj \vert P_jf(x_0)\vert < \infty.
$$
Hence, we can apply Lemma 3 to (3.2) for $t \ge T_1$, and we obtain
that $u_{\mu_1}(x_0,t) = u_{\mu_2}(x_0,t)$ for $T_1<t<T_2$ if and only if
$$
\int^T_0 e^{\la_ns}\mu(s) ds = 0, \quad n\in \N \setminus \Lambda(x_0).
                                                \eqno{(3.3)}
$$
We set $\eta = e^s$ and $h(\eta) = \mu(\log \eta)$, and we see that
$1 < \eta < e^T$ for $0<s<T$.  Then (3.3) is equivalent to
$$
\int^{e^T}_1 \eta^{\la_n}\frac{h(\eta)}{\eta} d\eta = 0,
\quad n\in \N \setminus \Lambda(x_0).                              \eqno{(3.4)}
$$
By the M\"untz theorem (e.g., E.7 on p.184 of
Borwein and Erd\'elyi \cite{BoE}) as directly applicable version, we see that
(3.4) yields $\frac{h(\eta)}{\eta} = 0$ for $ 1<\eta<e^T$ if and only if
(1.12) holds.  Thus the proof of Theorem 2 is complete.

In the case of $d=1$, we see that $\la_j = \sigma_j$ because all the
eigenvalues are simple.  Therefore Corollary follows directly
since $\la_j \sim \rho_0j^2$ as $j \to \infty$, which implies $\sumj \frac{1}{\la_j}
< \infty$ if $d=1$.
\section{Proofs of Theorem 3 and Proposition 1}

For general dimensions $d$, we have
$$
u_{\mu}(x,t) = \sumn \left(\int^t_0 \frac{\sin \sqrt{\la_n}(t-s)}{\sqrt{\la_n}}
\mu(s) ds\right) P_nf(x), \quad x \in \OOO, \, t>0
                                                    \eqno{(4.1)}
$$
(e.g., Komornik \cite{Ko}).  Setting $\mu := \mu_1 - \mu_2$ and
$u := u_{\mu_1} - u_{\mu_2}$, since $\mu(t) = 0$ for $t>T$ and
$u \in C(\ooo{\OOO} \times [0,\infty))$, in terms of (1.8) and
Lemma 1, we have
\begin{align*}
&u_{\mu}(x_0,t) = \sumn \left( \int^T_0 \frac{\sin\sqrt{\la_n}(t-s)}
{\sqrt{\la_n}}\mu(s) ds \right) P_nf(x_0)\\
= &\sumn \frac{\sin \sqrt{\la_n}t}{\sqrt{\la_n}}
\left( \int^T_0 \mu(s) \cos\sqrt{\la_n}s ds \right) P_nf(x_0)\\
- &\sumn \frac{\cos\sqrt{\la_n}t}{\sqrt{\la_n}}
\left( \int^T_0 \mu(s)\sin\sqrt{\la_n}s ds \right) P_nf(x_0),
\quad t>T.
\end{align*}
Therefore,
$$
u_{\mu}(x_0,t) = \sumn (a_n(P_nf)(x_0)\sin \sqrt{\la_n}t
+ b_n(P_nf)(x_0)\cos \sqrt{\la_n}t), \quad t>T,   \eqno{(4.2)}
$$
where we set
$$
a_n = \frac{1}{\sqrt{\la_n}}\int^T_0 \mu(s)\cos \sqrt{\la_n}s ds, \quad
b_n = \frac{-1}{\sqrt{\la_n}}\int^T_0 \mu(s)\sin \sqrt{\la_n}s ds, \quad
n\in N.
$$
The series is convergent in $C(\ooo{\OOO} \times [0,\infty)])$ by
(1.8) and Lemma 1 (i).
\\
{\bf Proof of Theorem 3.}

In the case of $\OOO = (0,\ell)$, we know
$d_n=1$, $\la_n = \frac{n^2\pi^2}{\ell^2}$ and
$\va_n(x) = \sqrt{\frac{2}{\ell}}\sin \frac{n\pi}{\ell}x$,
$(P_nf)(x_0) = (f, \va_n)\va_n(x_0)$.   Therefore,
$$
u_{\mu}(x_0,t) = \sumn \left( a_n\sin \frac{n\pi t}{\ell}
+ b_n \cos\frac{n\pi s}{\ell} ds \right)
(f,\va_n)\va_n(x_0) \quad \mbox{for $T_1<t<T_2$}
                                        \eqno{(4.3)}
$$
for $T_1<t<T_2$.

Since $T_2 - T_1 \ge 2\ell$ by (1.15),
we have (4.3) for $T_1<t<T_1+2\ell$.
Since (4.3) is convergent in $C[0,\infty)$, we see that
(4.3) is convergent in $L^2(T_1,T_1+2\ell)$.  Hence, taking
the scalar products in $L^2(T_1,T_1+2\ell)$ with $\sin\frac{m\pi t}{\ell}$
and $\cos\frac{m\pi t}{\ell}$, we obtain $a_m(f,\va_m)\va_m(x_0)
= b_m(f,\va_m)\va_m(x_0) = 0$ for $m\in\N$.
By (1.15), we see that $(f,\va_m) \ne 0$ and
$\va_m(x_0) = \sqrt{\frac{2}{\ell}} \sin \frac{m\pi}{\ell}x_0
\ne 0$ for all $m\in \N$.  Consequently, $a_m=b_m=0$, that is,
$$
\int^T_0 \mu(s)\cos\frac{m\pi s}{\ell} ds
= \int^T_0 \mu(s)\sin\frac{m\pi s}{\ell} ds = 0, \quad m\in \N.
$$
We recall the assumption $T \le 2\ell$.  If $T < 2\ell$, then we set
$\www{\mu}(t) = \left\{\begin{array}{rl}
\mu(s), \quad & 0<s<T,\\
0,      \quad & T\le s < 2\ell.
\end{array}\right.$
Then
$$
\int^{2\ell}_0 \www{\mu}(s)\cos\frac{m\pi s}{\ell} ds
= \int^{2\ell}_0 \www{\mu}(s)\sin\frac{m\pi s}{\ell} ds = 0, \quad m\in \N.
$$
Therefore, $\www{\mu}(s) = c_1$: constant for $0<s<2\ell$.
In the case of $T=2\ell$, we have $c_1 = \mu(T) = 0$ by
$\mu(t) = 0$ for $t \ge T$.  In the case of $T<2\ell$, since $\www{\mu}
= 0$ in $[T, 2\ell]$, we see $c_1=0$.  Hence $\mu(s) = 0$ for $0<s<T$, and
the proof of Theorem 3 is complete.
\\
\vspace{0.2cm}
{\bf Proof of Proposition 1.}\\
Since $f \in C^{\infty}_0(\OOO)$, we see (4.2), and the series
$u_{\mu}(x_0,t)$ is convergent in $L^{\infty}(\R)$, and so $\gamma$ is an almost
periodic function (e.g., B\"ottcher and Silbermann \cite{BoSi}).
By a property of almost periodic functions (e.g., property (d) in p.493
in \cite{BoSi}) yields
$$
\Vert u_{\mu}(x_0,\cdot)\Vert_{L^{\infty}(\R)}
= \limsup_{t\to\infty} \vert u_{\mu}(x_0,t)\vert.
$$
By (1.18), we see that $\Vert u_{\mu}(x_0,\cdot)\Vert_{L^{\infty}(\R)} = 0$,
that is,
$$
\sumn a_n(P_nf)(x_0)\sin \sqrt{\la_n}t
+ \sumn b_n(P_nf)(x_0)\cos \sqrt{\la_n}t = 0, \quad t\in \R.
$$
Substituting $-t$ and adding and subtracting, we obtain
$$
\sumn a_n(P_nf)(x_0)\sin \sqrt{\la_n}t
= \sumn b_n(P_nf)(x_0)\cos \sqrt{\la_n}t = 0, \quad t\in \R.
                                                              \eqno{(4.4)}
$$

Now we can prove\\
{\bf Lemma 4.}\\
{\it
Let $0<p_1<p_2< \cdots \, \longrightarrow \infty$.
%and the set
%$\{ p_n\}_{n\in \N}$ has no accumulation points except for $\infty$.
Let $\alpha_n \in \R$ satisfy
$$
\sumn \vert \alpha_n \vert < \infty.                \eqno{(4.5)}
$$
Then
$$
\sumn \alpha_n \sin p_nt = 0, \quad t>0 \quad \mbox{or}\quad
\sumn \alpha_n \cos p_nt = 0, \quad t>0
$$
yields $\alpha_n=0$ for each $n\in \N$.
}

For completeness the proof of the lemma is provided at the end of this section.
\\

Now we complete the proof of Proposition 1.
Applying Lemma 4 to (4.4), we obtain
$$
a_n(P_nf)(x_0) = b_n(P_nf)(x_0) = 0, \quad n\in \N.
$$
By the definition (1.10) of $\Lambda(x_0)$, we reach
$a_n=b_n=0$ for $n\in \N \setminus \Lambda(x_0)$, that is,
$$
\int^T_0 \mu(s) \cos\sqrt{\la_n}s ds = \int^T_0 \mu(s) \sin\sqrt{\la_n}s ds
= 0, \quad n \in \N \setminus \Lambda(x_0).
$$
With the even extension of $\mu$ to $(-T,0)$ for example, we see
$$
\int^T_{-T} \mu(s) e^{s\sqrt{-1}\sqrt{\la_n}}ds = 0 \quad
\mbox{for $n \in \N \setminus \Lambda(x_0)$}.
                                                     \eqno{(4.6)}
$$

We can prove the following lemma, whose proof is given at the end of
this section.
\\
{\bf Lemma 5.}\\
{\it
Under assumptions (1.16) and (1.17), the system
$$
\{ e^{s\sqrt{-1}\sqrt{\la_n}} \}_{n\in \N \setminus \Lambda(x_0)}
$$
is complete in $C[-T,T]$.  In particular, (4.6) implies
$\mu(s) = 0$ for $0<t<T$.
}

Here we say that
$\{ e^{s\sqrt{-1}\sqrt{\la_n}} \}_{n\in \N \setminus \Lambda(x_0)}$ is
complete in $C[-T,T]$ if $g \in C[-T,T]$ and
$$
\int^T_{-T} g(s) e^{s\sqrt{-1}\sqrt{\la_n}}ds = 0
\quad \mbox{for all $n\in \N \setminus \Lambda(x_0)$},
$$
imply $g(s) = 0$ for $0\le s\le T$.

By applying Lemma 5 to (4.6), the proof of Proposition 1 is finished.
\\

We conclude this section with
\\
{\bf Proof of Lemma 4.}
It suffices to prove in the case $\sumn \alpha_n\sin p_nt = 0$ for
$t>0$.  By (4.5), the series
$\sumn \alpha_n\sin p_nt = \frac{1}{2}\sumn \alpha_n(e^{\sqrt{-1}p_nt}
- e^{-\sqrt{-1}p_nt})$ is convergent in $L^{\infty}(\R)$.
Therefore, we can take the Laplace transform term-wisely:
$$
0 = \sumn \alpha_n \int^{\infty}_0 (e^{\sqrt{-1}p_nt}
- e^{-\sqrt{-1}p_nt})e^{-\xi t} dt
= \sumn \frac{2\sqrt{-1}\alpha_np_n}{\xi^2+p_n^2}
$$
for Re $\xi > 0$.  Setting $\eta = \xi^2$, we have
$$
\sumn \frac{\alpha_np_n}{\eta+p_n^2} = 0, \quad \eta \ge 0.
$$
For $\eta > 0$ and $k\in \N$, since $\vert \eta + p_n^2\vert^k
\ge p_n^{2k}$ and so
$$
\left\vert \frac{\alpha_np_n}{(\eta + p_n^2)^k}\right\vert
\le \frac{\vert \alpha_n\vert}{p_n^{2k-1}} \le C\vert \alpha_n\vert,
\quad n\in \N.
$$
In view of (4.5), the series $\sumn \frac{\alpha_np_n}{(\eta+p_n^2)^k}$
is uniformly convergent in $\eta \in [0,L]$ with arbitrary $L>0$ for
each $k\in \N$.
Therefore, we can term-wisely differentiate in $\eta$:
$$
\sumn \frac{\alpha_np_n}{(\eta + p_n^2)^k} = 0, \quad \eta > 0.
$$
Substituting $\eta = 1$, we can write
$$
\frac{\alpha_1p_1}{(1 + p_1^2)^k} + \sum_{k=2}^{\infty}
\frac{\alpha_np_n}{(1 + p_n^2)^k} = 0,
$$
that is,
$$
\frac{\alpha_1p_1}{1 + p_1^2}
= -\sum_{k=2}^{\infty} \alpha_n\frac{p_n}{1+p_n^2}
\left( \frac{1+p_1^2}{1+p_n^2}\right)^{k-1}.
$$
Since $0<p_1<p_2< \cdots$, we have
$$
\frac{1+p_1^2}{1+p_n^2} \le \frac{1+p_1^2}{1+p_2^2} =: r < 1
$$
for each $n\ge 2$ and $\frac{p_n}{1+p_n^2}\le 1$.  Hence we see that
$$
\left\vert \frac{\alpha_1p_1}{1 + p_1^2} \right\vert \le \sum_{n=2}^{\infty}
\vert \alpha_n\vert r^{k-1} \le Cr^{k-1}\quad \longrightarrow \quad 0
$$
as $k\to \infty$ in view of (4.5).
Therefore, $\alpha_1=0$ by $p_1\ne 0$.  Continuing this argument, we can
successively obtain $\alpha_n=0$ for $n\in \N$, so that the proof of Lemma 4
is complete. $\blacksquare$
\\
{\bf Proof of Lemma 5.}\\
Considering the zero extension, we readily see that
$\{ e^{s\sqrt{-1}\sqrt{\la_n}} \}_{n\in \N \setminus \Lambda(x_0)}$ is
complete in $C[-T,T]$, then so is it in $C[-\www{T},\, \www{T}]$ with
any $0<\www{T} < T$.

We set
$$
R:= \sup\{ L;\,
\{ e^{s\sqrt{-1}\sqrt{\la_n}} \}_{n\in \N \setminus \Lambda(x_0)}
\quad \mbox{is complete in $C[-L,L]$}\}
$$
and
$$
D:= \lim_{k\to \infty, \in \N \setminus \Lambda(x_0)}
\frac{k}{\sqrt{\la_k}}.
$$
Since $\sqrt{\la_k} > 0$ for $k \in \N \setminus \Lambda(x_0)$ and
is increasing in $k$, we can apply Theorem 13 (p.116) in \cite{Yo}, and we
can obtain
$$
\pi D \le R.
$$
Therefore, if $T < \pi D$, then we have $T < R$.  Hence,
$\{ e^{s\sqrt{-1}\sqrt{\la_n}} \}_{n\in \N \setminus \Lambda(x_0)}$
is complete in $C[-T,T]$ if $T < \pi D$.
We finished the proof of Lemma 5.
$\blacksquare$
\section{Proof of Theorem 4}
{\bf First Step:}
\\
We have
$$
\left\{ \begin{array}{rl}
& \ppp_tu(x,t) = \Delta u(x,t), \quad x\in \OOO, \, t>T,\\
& u(x,t) = 0, \quad x\in \ppp\OOO, \, t>T.
\end{array}\right.                       \eqno{(5.1)}
$$
We can prove
$$
\Vert u(\cdot,\www{t})\Vert_{L^2(\OOO)}
\le C\Vert \ppp_{\nu}u\Vert_{L^2(\gamma \times (T_1,T_2))} \quad
\mbox{for $T_1 < \www{t} < T_2$}.                                  \eqno{(5.2)}
$$
\\
{\bf Proof of (5.2).}
\\
First we show a Carleman estimate.  For it,
let $d \in C^2(\ooo{\OOO})$ satisfy $d>0$ in $\OOO$,
$\vert \nabla d\vert > 0$ on $\ooo{\OOO}$ and
$\ppp_{\nu}d \le 0$ on $\ppp\OOO \setminus \Gamma$.
For the proof of such $d$, see e.g., Imanuvilov \cite{Ima}.

We set
$$
\beta(x,t) = \frac{ e^{\lambda d(x)}}{(t-T_1)(T_2-t)}, \quad
\alpha(x,t) = \frac{ e^{\lambda d(x)}
- e^{2\lambda\Vert d\Vert_{C(\ooo{\OOO})}}}
{(t-T_1)(T_2-t)}, \quad x\in \OOO, \, T_1<t<T_2,
$$
where we choose $\lambda > 0$ sufficiently large.
Then
\\
{\bf Lemma 6 (Carleman estimate)}.\\
{\it
There exist constants $C>0$ and $s_0>0$ such that
$$
\int_{\OOO \times (T_1,T_2)} \left( \frac{1}{s\beta}
\vert \ppp_tu\vert^2 + s^3\la^4\beta^3\vert u\vert^2\right)
e^{2s\alpha} dxdt
\le C\int_{\gamma \times (T_1,T_2)} \vert \ppp_{\nu}u\vert^2 dS dt
$$
for all $s>s_0$ and all $u$ satisfying (5.1) and $\ppp_tu, \Delta u
\in L^2(\OOO \times (T_1,T_2))$.
}

For the proof, we refer to Imanuvilov \cite{Ima}, and also
Imanuvilov and Yamamoto \cite{IY1}, where the parabolic equation is considered
in the time interval $(0,T)$, and we can translate by change
$t \longmapsto \frac{T}{T_2-T_1}(t-T_1)$.
\\

We choose and fix $t_1, t_2$ such that
$T_1<t_1 < \frac{T_1+T_2}{2} < t_2 < T_2$.
Since $\alpha(x,t) < 0$ for $x \in \OOO$ and $T_1<t<T_2$,
there exist constants $C_1, C_2, C_3 > 0$ such that
$$
\alpha(x,t) \ge -C_1, \quad C_2 \le \beta(x,t)<C_3, \quad
(x,t) \in \ooo{\OOO} \times [t_1,t_2].
$$
Therefore, for all large $s>0$, we obtain
$$
\int_{\OOO \times (t_1,t_2)} (\vert \ppp_tu\vert^2
+ \vert u\vert^2) dxdt
\le C\int_{\gamma \times (T_1,T_2)} \vert \ppp_{\nu}u\vert^2
dS dt.
$$
Here the constant $C>0$ depends on $t_1, t_2$ and fixed large
constants $s>0$.
Applying the Sobolev embedding $H^1(t_1,t_2;L^2(\OOO)) \subset
C([t_1,t_2];L^2(\OOO))$, we obtain
$$
\sup_{t_1\le t\le t_2} \Vert u(\cdot,t)\Vert^2_{L^2(\OOO)}
\le C\Vert \ppp_{\nu}u\Vert^2_{L^2(\gamma \times (t_1,t_2))}.
$$
Since $t_1, t_2$ are atbitrarily provided that
$T_1 < t_1 < t_2 < T_2$, the proof of (5.2) is complete.
\\

Henceforth we renumber the set $\{ \la_k\}_{k\in \N}$
of all the eigenvalues of $-\Delta$ with the zero
Dirichlet boundary condition, according to the multiplicities:
$0<\sigma_1 \le \sigma_2 \le \sigma_3 \le \cdots$.
We can choose an eigenfunction $\va_k$ corresponding to $\sigma_k$
such that $\{ \va_k\}_{k\in \N}$ forms an orthonormal basis in $L^2(\OOO)$.

By the eigenfunction expansion, we can represent the solution
$u$ to (1.1) for $0\le t \le T$ as
$$
u(x,t) = \sumj e^{-\sigma_jt}\left(\int^t_0 e^{\sigma_js}\mu(s) ds\right)
(f, \va_j)\va_j(x), \quad x \in \OOO, \, 0\le t \le T.      \eqno{(5.3)}
$$
Since $\mu(s) = 0$ for $s \ge T$, we obtain
$$
u(x,t) = \sumj e^{-\sigma_jt}\left(\int^T_0 e^{\sigma_js}\mu(s) ds\right)
(f, \va_j)\va_j(x), \quad x \in \OOO, \, T<t\le T_2.      \eqno{(5.4)}
$$
\\
{\bf Second Step: Proof of Theorem 4 (i).}
\\
By (5.2) and $\ppp_{\nu}u=0$ on $\gamma \times (T_1,T_2)$, we have
$u(x,t) = 0$ for $x\in \OOO$ and $T_1 < t < T_2$.  Therefore, (5.4) yields
$$
\sumj e^{-\sigma_j\www{t}}\left(\int^T_0 e^{\sigma_js}\mu(s) ds\right)
(f, \va_j)\va_j(x) = 0, \quad x \in \OOO,\, T_1 \le \www{t} \le T_2,
$$
which implies
$$
\left(\int^T_0 e^{\sigma_js}\mu(s) ds\right)
(f, \va_j)= 0 \quad \mbox{for all $j\in \N$}.
$$
By the assumption (1.29), we reach $(f, \va_j) = 0$ for all
$j\in \N$.  Thus $f=0$ in $\OOO$.
\\
{\bf Third Step: Proof of Theorem 4 (ii).}
\\
We fix $\www{t} \in (T_1,T_2)$ arbitrarily.
Combining (5.2) with (5.4) with $t=\www{t}$, we obtain
$$
\Vert u(\cdot,\www{t})\Vert_{L^2(\OOO)}^2
= \left\Vert \sumj p_j(f,\va_j)\va_j\right\Vert^2_{L^2(\OOO)}
\le C\Vert \ppp_{\nu}u\Vert^2_{L^2(T_1,T_2;L^2(\gamma))}.
                                                       \eqno{(5.5)}
$$
Here we set
$$
p_j = e^{-\sigma_j\www{t}}\left(\int^T_0 e^{\sigma_js}\mu(s) ds\right),
\quad j\in \N.
$$
Since $f\in \mathcal{D}(A^{\ell})$ with given $\ell \in \N$, we have
$$
\vert (f,\va_j)\vert = \vert (A^{\ell}f, A^{-\ell}\va_j)\vert
= \left\vert \frac{1}{\sigma_j^{\ell}}(A^{\ell}f,\va_j)\right\vert.
%\le \frac{1}{\sigma_j^{\ell}}\Vert A^{\ell}f\Vert_{L^2(\OOO)}.
                                           \eqno{(5.6)}
$$
Therefore, for each $N\in \N$, we obtain
\begin{align*}
&\sumj \vert (f,\va_j)\vert^2
= \left( \sum_{j=1}^N + \sum_{j=N+1}^{\infty}\right)
\vert (f, \va_j)\vert^2
= \sum_{j=1}^N \vert (f,\va_j)\vert^2
+ \sum_{j=N+1}^{\infty} \frac{1}{\sigma_j^{2\ell}}
\vert (A^{\ell}f, \va_j)\vert^2\\
\le& \sum_{j=1}^N \vert (f,\va_j)\vert^2
+ \frac{1}{\sigma_N^{2\ell}}\sum_{j=N+1}^{\infty}
\vert (A^{\ell}f, \va_j)\vert^2
\end{align*}
$$
\le \sum_{j=1}^N \vert (f,\va_j)\vert^2
+ \frac{1}{\sigma_N^{2\ell}}\Vert A^{\ell}f\Vert^2_{L^2(\OOO)}.
                                      \eqno{(5.7)}
$$

Next we estimate $p_j$.  By (1.30), we arbitrarily fix constants
$0<t_1<t_2<T$ and $\delta>0$ such that $\mu(s) \ge \delta$ for
$t_1\le s \le t_2$.  Therefore,
$$
\int^T_0 e^{\sigma_js}\mu(s) ds \ge \delta\int^{t_2}_{t_1} e^{\sigma_js} ds
= \frac{e^{\sigma_jt_2} - e^{\sigma_jt_1}}{\sigma_j}\delta,
$$
so that
$$
\frac{1}{\vert p_j\vert}
= \frac{e^{\sigma_j\www{t}}\sigma_j}{\delta}
\frac{1}{e^{\sigma_jt_2}(1-e^{-\sigma_j(t_2-t_1)})}
= \frac{e^{\sigma_j(\www{t}-t_2)}\sigma_j}{\delta}
\frac{1}{1 - e^{-\sigma_1(t_2-t_1)}} \le C_1\sigma_j e^{C_2\sigma_j}
                                                \eqno{(5.8)}
$$
for $j\in \N$.
By $e^{-\sigma_j(t_2-t_1)} \le e^{-\sigma_1(t_2-t_1)}
< 1$ for all $j\in \N$,
the constants $C_1>0$ and $C_2>0$ are independent of $j\in \N$.
On the other hand, substituting $t = \www{t}$ in (5.4),
we take the scalar product of (5.4) with $\va_j$, and we obtain
$(u(\cdot,\www{t}), \va_j) = p_j(f,\va_j)$.
Hence, (5.8) implies
$$
\vert (f, \va_j)\vert = \frac{1}{\vert p_j\vert}
\vert (u(\cdot,\www{t}), \va_j)\vert
\le C_1\sigma_j e^{C_2\sigma_j} \vert (u(\cdot,\www{t}), \va_j)\vert,
\quad j\in \N.
$$
Substituting this into (5.7), we see
\begin{align*}
& \Vert f\Vert^2_{L^2(\OOO)}
\le C_1^2 \sum_{j=1}^N \sigma_j^2 e^{2C_2\sigma_j}
\vert (u(\cdot,\www{t}), \va_j)\vert^2
+ \frac{1}{\sigma_N^{2\ell}}\Vert A^{\ell}f\Vert^2_{L^2(\OOO)}\\
\le& C_1^2 \sigma_N^2 e^{2C_2\sigma_N}\Vert (u(\cdot,\www{t})\Vert^2
_{L^2(\OOO)}
+ \frac{1}{\sigma_N^{2\ell}}\Vert A^{\ell}f\Vert^2_{L^2(\OOO)}
\end{align*}
with arbitrary $N\in \N$.

In terms of (5.2), since we are considering the case of
$\Vert \ppp_{\nu}u\Vert_{L^2(\gamma \times (T_1,T_2))} \longrightarrow 0$,
we can assume that
$\Vert u(\cdot,\www{t})\Vert_{L^2(\OOO)} < 1$.
Since $\sigma_N \sim \rho_0N^{\frac{2}{d}}$ as $N\to \infty$
(e.g., \cite{Ag}, \cite{CH}) and $N^{\frac{2}{d}}
\le \exp(N^{\frac{2}{d}})$ for sufficiently large $N\in \N$, we obtain
$$
\Vert f\Vert_{L^2(\OOO)} \le C_1\exp\left( C_3N^{\frac{2}{d}}\right)\eta
+ \frac{C_1}{N^{\theta_0}}M \quad \mbox{for all $N\in \N$},
                                                             \eqno{(5.9)}
$$
where we choose a small constant $\eta_0\in (0,1)$ and set
$$
\eta := \Vert u(\cdot,\www{t})\Vert_{L^2(\OOO)} \le 1 - \eta_0, \quad
M := \Vert A^{\ell}f\Vert_{L^2(\OOO)}, \quad
\theta_0 := \frac{2\ell}{d}.
$$
Now, for given $\eta > 0$,
we make the right-hand side of (5.9) smaller by choosing $N\in \N$.
By $\eta < 1$, we see $\log\frac{1}{\eta} > 1$.
Let $\theta \in (0, \, \ell)$ be arbitrarily given.
We can choose $N = N(\eta) \in \N$ such that
$$
\left( \log\frac{1}{\eta}\right)^{\frac{\theta}{\theta_0}}
\le N < \left( \log\frac{1}{\eta}\right)^{\frac{\theta}{\theta_0}}
+ 1.
$$
Then
$$
\left( \log\frac{1}{\eta}\right)^{\theta} \le N^{\theta_0},
$$
that is,
$$
\frac{C_1}{N^{\theta_0}} \le \frac{C_1}
{\left( \log\frac{1}{\eta}\right)^{\theta}}.                 \eqno{(5.10)}
$$
Moreover, we know
$$
N^{\frac{2}{d}} \le \left( \left( \log\frac{1}{\eta}\right)
^{\frac{\theta}{\theta_0}} + 1 \right)^{\frac{2}{d}}
\le C_4\left( \left( \log\frac{1}{\eta}\right)
^{\frac{\theta}{\theta_0}\frac{2}{d}} + 1 \right),
$$
that is,
$$
e^{C_3N^{\frac{2}{d}}}\eta^{\frac{1}{2}}
\le \eta^{\frac{1}{2}}\exp\left( C_3C_4\left( \left(
\log\frac{1}{\eta}\right)^{\frac{\theta}{\ell}} + 1 \right)\right)
\le C_5\eta^{\frac{1}{2}} \exp\left( C_5
\left(\log\frac{1}{\eta}\right)^{\frac{\theta}{\ell}}\right)
$$
for $\eta \le 1 - \eta_0$.
Setting $t = \log\frac{1}{\eta} > 0$, we can see
$$
\eta^{\frac{1}{2}} \exp\left( C_5
\left(\log\frac{1}{\eta}\right)^{\frac{\theta}{\ell}}\right)
= e^{-\frac{1}{2}t}\exp(C_5t^{\frac{\theta}{\ell}}),
$$
and so $\frac{\theta}{\ell} < 1$ yields
$$
\sup_{\eta<1} \left\vert \eta^{\frac{1}{2}} \exp\left( C_5
\left(\log\frac{1}{\eta}\right)^{\frac{\theta}{\ell}}\right)\right\vert
< \infty.
$$
Hence,
$$
\exp(C_3N^{\frac{2}{d}})\eta
\le C_5\eta^{\frac{1}{2}}
\sup_{\eta < 1} \left\vert \eta^{\frac{1}{2}} \exp\left( C_5
\left(\log\frac{1}{\eta}\right)^{\frac{\theta}{\ell}}\right)\right\vert
\le C_6\eta^{\frac{1}{2}},
$$
and (5.9) with this choice of $N$ and (5.10) yield
$$
\Vert f\Vert_{L^2(\OOO)} \le C_7\eta^{\frac{1}{2}}
+ \frac{C_1}{\left(\log\frac{1}{\eta}\right)^{\theta}}M.
$$
Since we can find a constant $C_8>0$ such that
$\eta^{\frac{1}{2}} \le \frac{C_8}
{\left(\log\frac{1}{\eta}\right)^{\theta}}$ for $\eta < 1$,
we obtain
$$
\Vert f\Vert_{L^2(\OOO)} \le \frac{C_9}
{\left( \log\frac{1}{\eta} \right)^{\theta}}M.
$$
By using
(5.5), we complete the proof of Theorem 4 (ii).
\section{Proof of Theorem 5}

Since $\mu(t) = 0$ for $t > T$, we have
$$
\left\{ \begin{array}{rl}
& \ppp_t^2u(x,t) = \Delta u(x,t), \quad x\in \OOO, \, t>T,\\
& u(x,t) = 0, \quad x\in \ppp\OOO, \, t>T.
\end{array}\right.                       \eqno{(6.1)}
$$
Setting $v = \ppp_tu$ and
$$
\www{\mu}(t) =
\left\{ \begin{array}{rl}
\mu'(t), \quad & 0 < t \le T, \\
0, \quad & t>T,
\end{array}\right.
$$
we have
$$
\left\{ \begin{array}{rl}
& \ppp_t^2v(x,t) = \Delta v(x,t) + \www{\mu}(t)f(x), \quad x\in \OOO, \, t>0,\\
& v(x,t) = 0, \quad x\in \ppp\OOO, \, t>0, \\
& v(x,0) = \ppp_tv(x,0) = 0, \quad x \in \OOO.
\end{array}\right.                                \eqno{(6.2)}
$$
By Theorem 8.2 (p.275) in \cite{LM}, we note that
$v\in C([0,\infty);H^1_0(\OOO))$ and $\ppp_tv \in
C([0,\infty);L^2(\OOO))$.  In particular,
$$
\left\{ \begin{array}{rl}
& \ppp_t^2v(x,t) = \Delta v(x,t), \quad x\in \OOO, \, t>T,\\
& v(x,t) = 0, \quad x\in \ppp\OOO, \, t>0, \\
& v(x,T) = \ppp_tu(x,T), \quad \ppp_tv(x,0) = \Delta u(x,T),
\quad x \in \OOO.
\end{array}\right.                  \eqno{(6.3)}
$$
In view of (1.24), we can apply the observability inequality
(e.g., Komornik \cite{Ko}) to (6.3),
so that we obtain
\begin{align*}
& \Vert \ppp_tu(\cdot,T)\Vert_{H^1(\OOO)} + \Vert \Delta u(\cdot,T)\Vert
_{L^2(\OOO)}
= \Vert v(\cdot,T)\Vert_{H^1(\OOO)} + \Vert \ppp_tv(\cdot,T)\Vert
_{L^2(\OOO)}\\
\le& C\Vert \ppp_{\nu}u\Vert_{H^1(T_1,T_2;L^2(\Gamma))}.
\end{align*}
Since $u(\cdot,T) = 0$ on $\ppp\OOO$, the elliptic regularity implies
$\Vert u(\cdot,T)\Vert_{H^2(\OOO)} \le C\Vert \Delta u(\cdot,T)\Vert
_{L^2(\OOO)}$, we have
$$
\Vert u(\cdot,T)\Vert_{H^2(\OOO)} + \Vert \ppp_tu(\cdot,T)\Vert
_{H^1_0(\OOO)} \le C\Vert \ppp_{\nu}u\Vert_{H^1(T_1,T_2;L^2(\Gamma))}.
                                                 \eqno{(6.4)}
$$
Similarly to (4.1), we can obtain
$$
u(x,t) = \sumn \left(\int^t_0 \frac{\sin\sqrt{\la_n}(t-s)}{\sqrt{\la_n}}
\mu(s) ds\right) P_nf(x)
$$
in $H^2(\OOO) \cap H^1_0(\OOO)$ for $0<t<T$.
Therefore,
$$
\ppp_tu(x,t) = \sumn \left(\int^t_0 \cos\sqrt{\la_n}(t-s)\mu(s) ds\right)
P_nf(x)
$$
in $H^1_0(\OOO)$ for $0<t<T$.
Hence,
$$
\left\{ \begin{array}{rl}
& u(x,T) = \sumn \left(\int^T_0 \frac{\sin\sqrt{\la_n}(T-s)}{\sqrt{\la_n}}
\mu(s) ds\right) P_nf(x), \\
& \ppp_tu(x,T) = \sumn \left(\int^T_0 \cos\sqrt{\la_n}(T-s)\mu(s) ds\right)
P_nf(x) \quad \mbox{in $L^2(\OOO)$}.
\end{array}\right.
                                               \eqno{(6.5)}
$$
On the other hand, by $\mu(T) = 0$, the integration by parts yields
$$
\int^T_0 \mu(s)\sin\sqrt{\la_n}(T-s) ds
= -\frac{\mu(0)\cos T\sqrt{\la_n}}{\sqrt{\la_n}}
- \frac{1}{\sqrt{\la_n}} \int^T_0 \mu'(s)\cos\sqrt{\la_n}(T-s) ds.
$$
The Riemann-Lebesgue theorem implies
$$
\int^T_0 \mu(s) \sin\sqrt{\la_n}(T-s) ds
= -\frac{\mu(0)\cos T\sqrt{\la_n}}{\sqrt{\la_n}}
+ o\left( \frac{1}{\sqrt{\la_n}} \right) \quad \mbox{as $n\to \infty$}.
                                                        \eqno{(6.6)}
$$
Similarly we can verify
$$
\int^T_0 \mu(s)\cos\sqrt{\la_n}(T-s) ds
= \frac{\mu(0)\sin T\sqrt{\la_n}}{\sqrt{\la_n}}
+ o\left( \frac{1}{\sqrt{\la_n}} \right) \quad \mbox{as $n\to \infty$}.
                                                        \eqno{(6.7)}
$$
Moreover, by (6.5), we have
\begin{align*}
& \Vert u(\cdot,T)\Vert^2_{H^2(\OOO)} + \Vert \ppp_tu(\cdot,T)\Vert^2
_{H^1_0(\OOO)}\\
=& \sumn \la_n\left( \left\vert
\int^T_0 \mu(s)\sin\sqrt{\la_n}(T-s) ds\right\vert^2
+ \left\vert \int^T_0 \mu(s)\cos\sqrt{\la_n}(T-s) ds\right\vert^2\right)
\Vert P_nf\Vert^2_{L^2(\OOO)}.
\end{align*}
By (6.6) and (6.7), we see
\begin{align*}
J_n &:= \la_n\left( \left\vert
\int^T_0 \mu(s)\sin\sqrt{\la_n}(T-s) ds\right\vert^2
+ \left\vert \int^T_0 \mu(s)\cos\sqrt{\la_n}(T-s) ds\right\vert^2\right)\\
=& \mu(0)^2(\cos^2 T\sqrt{\la_n} + \sin^2 T\sqrt{\la_n})
+ o(1) \ge \mu(0)^2 - o(1)
\end{align*}
as $n \longrightarrow \infty$.
We choose large $N\in \N$ such that $J_n \ge \frac{\mu(0)^2}{2}$
for $n \ge N$.  By (1.31) we obtain
$$
\mu_0:= \min\left\{ \min_{1\le n\le N}J_n, \, \frac{\mu(0)^2}{2} \right\}
> 0,
$$
and so
$$
\la_n\left( \left\vert \int^T_0 \mu(s)\sin\sqrt{\la_n}(T-s) ds\right\vert^2
+ \left\vert \int^T_0 \mu(s)\cos\sqrt{\la_n}(T-s) ds\right\vert^2\right)
\ge \mu_0, \quad n \in \N,
$$
which yields
$$
\Vert u(\cdot,T)\Vert^2_{H^2(\OOO)}
+ \Vert \ppp_tu(\cdot,T)\Vert^2_{H^1_0(\OOO)}
\ge \mu_0 \sumn \sum_{k=1}^{d_n} \vert (f, \va_{nk})\vert^2
= \mu_0\Vert f\Vert^2_{L^2(\OOO)}.
$$
The combination with (6.4) completes the proof of the second inequality in
(1.32).

Finally we prove the first inequality in (1.32).  Applying the direct
inequality (e.g., Komornik \cite{Ko}) to (6.1), for any
$T_0 > T$, we have
$$
\Vert \ppp_{\nu}u\Vert_{H^1(T,T_0;L^2(\ppp\OOO))}
\le \Vert \ppp_{\nu}v\Vert_{L^2(T,T_0;L^2(\ppp\OOO))}
\le C(T_0)(\Vert v(\cdot,T)\Vert_{H^1(\OOO)}
+ \Vert \ppp_tv(\cdot,T)\Vert_{L^2(\OOO)}).
$$
Here $C(T_0)>0$ depends on $T_0 > T$.  Hence,
$$
\Vert \ppp_{\nu}u\Vert_{H^1(T_1,T_2;L^2(\ppp\OOO))}
\le C(\Vert v(\cdot,T)\Vert_{H^1(\OOO)}
+ \Vert \ppp_tv(\cdot,T)\Vert_{L^2(\OOO)}).             \eqno{(6.8)}
$$
By $\mu' \in L^2(0,T)$, we apply the usual energy estimate
(or Theorem 8.2 (p.275) in \cite{LM} for example) to (6.2) for
$0<t<T$, we have
$$
\Vert v(\cdot,T)\Vert_{H^1(\OOO)} + \Vert \ppp_tv(\cdot,T)\Vert_{L^2(\OOO)}
\le C\Vert \mu'f\Vert_{L^2(0,T;L^2(\OOO))},
$$
with which (6.8) completes the proof of the first inequality in (1.32).
Thus the proof of Theorem 5 is complete.
\section{Proofs of Propositions 2 and 3}

\subsection{Proof of Proposition 2}

We set $u:= u_{\mu_1} - u_{\mu_2}$ and $\mu:= \mu_1 - \mu_2$.
Then we have (5.1) for $t > T$.  Therefore, in view of (5.2), we obtain
$u\left( \cdot, \frac{T_1+T_2}{2}\right) = 0$ in $\OOO$.
The uniqueness for a heat equation (5.1) backward in time (e.g.,
Imanuvilov and Yamamoto \cite{IY2}, Isakov \cite{Is})
yields $u(\cdot,T) = 0$ in $\OOO$.  Here, by
$u \in C([0,T_2];L^2(\OOO))$, we see that
$$
\lim_{\ep\to 0, \ep>0} u(\cdot, T+\ep)
= \lim_{\ep\to 0, \ep>0} u(\cdot, T-\ep),
$$
and (3.1) yields
$$
\sumn e^{-\la_nT}\left(\int^T_0 e^{\la_ns}\mu(s) ds\right) (P_nf)(x) = 0, \quad
x\in \OOO.
$$
Therefore, by definition (1.21) of $\www{\Lambda}$, we obtain
$$
\int^T_0 e^{\la_ns}\mu(s) ds = 0, \quad n \in \N \setminus \www{\Lambda}.
$$
Similarly to the argument after (3.3), we apply the M\"untz theorem
(\cite{BoE}) to
obtain $\mu(t) = 0$ for $0<t<T$ if and only if (1.22) holds.
Thus the proof of Proposition 2 is complete.

\subsection{Proof of Proposition 3}

We set $u:= u_{\mu_1} - u_{\mu_2}$ and $\mu:= \mu_1 - \mu_2$.
Then we have (6.1).  First, since $\ppp_{\nu}u=0$ on
$\Gamma \times (T_1,T_2)$, by (1.24) the observability inequality
(e.g., \cite{Ko}) to (6.1) yields
$u(\cdot,T_1) = \ppp_tu(\cdot,T_1) = 0$ in $\OOO$.
Therefore, in (6.1) we see that
$u(\cdot,T) = \ppp_tu(\cdot,T) = 0$ in $\OOO$.
Here we note that $u \in C^1([0,T_2];L^2(\OOO))$.  Hence, (6.4) implies
$$
\int^T_0 \mu(s) \sin\sqrt{\la_n}(T-s)  ds = 0, \quad
\int^T_0 \mu(s) \cos\sqrt{\la_n}(T-s) \mu(s) ds = 0, \quad
n\in \N \setminus \www{\Lambda},
$$
that is,
$$
\left\{ \begin{array}{rl}
& \sin T\sqrt{\la_n}\int^T_0 \mu(s)\cos\sqrt{\la_n}s ds
- \cos T\sqrt{\la_n}\int^T_0 \mu(s)\sin\sqrt{\la_n}s ds = 0, \\
& \cos T\sqrt{\la_n}\int^T_0 \mu(s)\cos\sqrt{\la_n}s ds
+ \sin T\sqrt{\la_n}\int^T_0 \mu(s) \sin\sqrt{\la_n}s ds = 0, \quad
n\in \N \setminus \www{\Lambda}.
\end{array}\right.
$$
Considering the system as linear equations in
$\int^T_0 \mu(s)\cos\sqrt{\la_n}s ds$ and
$\int^T_0 \mu(s)\sin\sqrt{\la_n}s ds$, since the determinant of
the coefficient matrix is non-zero:
$\sin^2 T\sqrt{\la_n} + \cos^2 T\sqrt{\la_n} = 1\ne 0$, we see
$$
\int^T_0 \mu(s)\cos\sqrt{\la_n}s ds = \int^T_0 \mu(s) \sin\sqrt{\la_n}s ds
= 0, \quad n\in\N \setminus \www{\Lambda} .                       \eqno{(7.1)}
$$
In view of (1.23) and (1.25), we can argue similarly to the proof of
Proposition 1, and by Lemma 5, we can reach
$\mu(s) = 0$, $0<s<T$.  Thus the proof of Proposition 3 is complete.
\section{Appendix I. Uniqueness in determining $\mu(t)$ by data
over time interval $(0,T)$}

Our main interest is the inverse source problems by data over
the time interval $(T_1, T_2)$ where $T<T_1<T_2$, because
the inverse problems with data over $(0,T)$ have been well studied.
However, uniqueness results by data over $(0,T)$ seem a kind of
folklore, and it is not easy to find relevant articles.
Thus for completeness, we here show the uniqueness limited to
the case of the determination of $\mu(t)$.

Our argument can work for the case where $\Delta$ is replaced by
a general uniformly elliptic operator with some regularity
condition on the coefficients, but we are restricted to
(1.1) and (1.2) with the zero Dirichlet boundary condition.
For simplicity we assume that $f \in C^{\infty}_0(\OOO)$.
\\
\vspace{0.2cm}
{\bf Case 1: one-dimensional case and pointwise
data $u(x_0,t)$, $0<t<T$
with $x_0 \in \OOO:= (0, \ell)$.}

We recall
$$
\la_n = \frac{n^2\pi^2}{\ell^2}, \quad
\va_n(x) = \frac{\sqrt{2}}{\sqrt{\ell}}\sin \frac{n\pi}{\ell}x,
\quad n \in \N,
$$
and $(f,g) := \int^{\ell}_0 f(x)g(x) dx$ for $f,g \in L^2(0,\ell)$.
By $u_k$, $k=1,2$, we denote the solution to (1.1) and (1.2) respectively.
that is, $u_1$ and $u_2$ are the solutions for the diffusion and the
wave equations respectively.
\\
{\bf Proposition A.}\\
{\it
We assume
$$
\mbox{there exists $n_1 \in \N$ such that $(f,\va_{n_1})\va_{n_1}(x_0)
\ne 0$}.                                       \eqno{(8.1)}
$$
{\bf (i): Case of diffusion equation.}
If $u_1(x_0,t) = 0$ for $0<t<T$, then $\mu(t) = 0$ for $0<t<T$.
\\
{\bf (ii): Case of wave equation.}
If $u_2(x_0,t) = 0$ for $0<t<T$, then $\mu(t) = 0$ for $0<t<T-\ell$.
}

We can compare Proposition A (i) and (ii) with Theorems 2 and 3
respectively, which naturally require more conditions.
In (ii), if $T \le \ell$, then we can not conclude the uniqueness for any
time interval $(0,\ep)$ with arbitrary $\ep>0$,
and in view of the propagation speed $1$, it is natural that
$T>\ell$ is needed for the meaningful uniqueness.
\\
{\bf Remark.}\\
If either
$$
f(x_0) \ne 0,
$$
or
$$
\frac{x_0}{\ell} \not\in \Q, \quad f\not\equiv 0 \quad \mbox{in
$(0,\ell)$},
$$
then (8.1) holds.
Indeed if (8.1) does not hold, then $(f,\va_n)\va_n(x_0) = 0$ for all
$n\in \N$, which implies $f(x_0) = \sumn (f,\va_n)\va_n(x_0) = 0$.
Hence, in the first case, (8.1) holds.  Furthermore,
in the second case we note that
$\frac{x_0}{\ell} \not\in \Q$ implies that $\va_n(x_0)\ne 0$
for all $n\in \N$.
\\
\vspace{0.2cm}
{\bf Case 2: general dimensional case and data $\ppp_{\nu}u$
on lateral subboundary of $\ppp\OOO \times (0,T)$.}

For arbitrarily chosen $x_0 \in \R^d$, we recall
that $\nu=\nu(x)$ denotes the unit outward normal vector to
$\ppp\OOO$ at $x$, and
$$
\Gamma(x_0) := \{ x\in \ppp\OOO;\,
(x-x_0)\cdot\nu(x) \ge 0 \}, \quad R(x_0) := \max_{x\in \ooo{\OOO}}
\vert x-x_0\vert.
$$
\\
{\bf Proposition B.}\\
{\it
We assume that $f\not\equiv 0$ in $\OOO$.
\\
{\bf (i): Case of diffusion equation.}
Let $\gamma \subset \ppp\OOO$ be an arbitrarily fixed subboundary.
If $\ppp_{\nu}u_1\vert_{\gamma \times (0,T)} = 0$, then
$\mu(t) = 0$ for $0<t<T$.
\\
{\bf (ii): Case of wave equation.}
If $\ppp_{\nu}u_2\vert_{\Gamma(x_0) \times (0,T)} = 0$, then
$\mu(t) = 0$ for $0 < t < T-R(x_0)$.
}
\\

Similarly to Proposition A, in the case of the wave equation,
owing to the finite propagation speed $1$,
the time interval for the uniqueness of $\mu(t)$ is
reduced by $R(x_0)$ from the observation time length $T$.
The conclusion of (ii) makes sense only if $T > R(x_0)$.
Proposition B (i) and (ii) should be compared with Propositions 2 and 3
respectively where data are taken over $(T_1, T_2)$ with $T_1>0$.
Also in (ii), for the uniqueness we can take arbitrary subboundary with
relevantly long time interval as long as the wave equation has analytic
coefficients.   For the uniqueness we can apply Fritz John's global
Holmgren theorem (e.g., John \cite{J}), but we omit the details.

For the proofs of the propositions, we reduce the problems to
initial boundary value problems for equations with the zero
right-hand sides.  More precisely, for $k=1,2$, let $v_k = v_k(x,t)$ satisfy
$$
\left\{ \begin{array}{rl}
& \ppp_t^kv_k(x,t) = \Delta v_k(x,t), \quad x\in \OOO, \, t>0,\\
& v_k(x,t) = 0, \quad x\in \ppp\OOO, \, t>0, \\
&
\left\{ \begin{array}{rl}
&v_k(x,0) = f(x), \quad x \in \OOO, \quad \mbox{if $k=1$}, \\
&v_k(x,0) = 0, \quad \ppp_tv_k(x,0) = f(x), \quad
x \in \OOO, \quad \mbox{if $k=2$}. \\
\end{array}\right.
\end{array}\right.                       \eqno{(8.2)}
$$
Then, as is easily verified, we see
$$
u_k(x,t) = \int^t_0 \mu(t-s)v_k(x,s) ds, \quad x \in \OOO, \, t>0
                                                   \eqno{(8.3)}
$$
for $k=1,2$.
\\
{\bf Proof of Proposition A.}
\\
For $k=1,2$, let $u_k(x_0,t) = 0$ for $0<t<T$.  Then (8.3) yields
$$
\int^t_0 \mu(t-s)v_k(x_0,s) ds = 0, \quad 0<t<T.
$$
The Titchmarsh theorem on the convolution (Titchmarsh \cite{T}) yields
that there exists $t_* \in (0,T)$ such that
$$
\left\{ \begin{array}{rl}
& \mu(s) = 0, \quad 0<s<t_*, \\
& v_k(x_0,s) = 0, \quad 0<s<T-t_*.
\end{array}\right.                          \eqno{(8.4)}
$$
We assume that $t_* < T$.  Otherwise $\mu(s) = 0$ for $0<s<T$ have been
already proved.
\\
{\bf Proof of (i).}\\
We know
$$
v_1(x_0,t) = \sumn e^{-\la_nt}(f,\va_n)\va_n(x_0), \quad t>0.
$$
Then
$$
\sumn e^{-\la_nt}(f,\va_n)\va_n(x_0) = 0, \quad 0<t<T-t_*.
                                                   \eqno{(8.5)}
$$
We can apply Lemma 3 in Section 3 to obtain
$(f,\va_n)\va_n(x_0) = 0$ for $n\in \N$.
This contradicts assumption (8.1).  Hence $t_* < T$ is impossible and so
$t_*=T$, that is, $\mu(t) = 0$ for $0<t<T$.  Thus the proof of
Proposition A (i) is complete.
\\
{\bf Proof of (ii).}\\
In the case of the wave equation, we have
$$
v_2(x,t) = \sumn \frac{\sin\sqrt{\la_n}t}{\sqrt{\la_n}}
(f,\va_n)\va_n(x_0) = 0                        \eqno{(8.6)}
$$
for $0<t<T-t_*$ (e.g., Komornik \cite{Ko}).
If $T-t_* \ge \ell$, then (8.6) holds for $0<t<\ell$.
Recalling that
$\va_n(x) = \frac{\sqrt{2}}{\sqrt{\ell}} \sin \frac{n\pi}{\ell}t$
with $n\in \N$, we can obtain
$$
(f,\va_n)\va_n(x_0) = 0 \quad \mbox{for all $n\in \N$},
$$
which contradicts assumption (8.1).
Hence, $T - t_* < \ell$, that is, $t_* > T-\ell$.  Thus
(8.4) yields $\mu(t) = 0$ for $0<t<T-\ell$, and the proof of
Proposition A (ii) is complete.
\\
{\bf Proof of Proposition B.}
\\
Similarly to the proof of Proposition A, we obtain
$$
\left\{ \begin{array}{rl}
& \mu(s) = 0, \quad 0<s<t_*, \\
& \ppp_{\nu}v_k = 0 \quad \mbox{on $\gamma \times (0,T-t_*)$ if $k=1$},\\
& \ppp_{\nu}v_k = 0 \quad \mbox{on $\Gamma(x_0) \times (0,T-t_*)$ if $k=2$}.
\end{array}\right.                          \eqno{(8.7)}
$$
\\
{\bf Proof of (i).}\\
Assume that $t_* < T$.  Then $\ppp_{\nu}v_1 = 0$ on
$\gamma \times (0,T-t_*)$.  With this, we apply the classical
unique continuation for the heat equation (e.g., Isakov \cite{Is},
Yamamoto \cite{Ya}), and so $v_1 = 0$ in $\OOO \times (0,T)$, which yields
$f=0$ in $\OOO$.  This contradicts the assumption $f\not\equiv 0$ in $\OOO$.
Therefore, $t_* < T$ is impossible, and we reach $t_*=T$.
The proof of (i) is complete.
\\
{\bf Proof of (ii).}\\
Assume that $t_*<T$.  Then $\ppp_{\nu}v_2 = 0$ on
$\Gamma(x_0) \times (0,T-t_*)$.  If $T-t_* > R(x_0)$, then the
observability inequality (e.g., \cite{Ko}) or the unique continuation by
Carleman estimate (e.g., \cite{BY}) yields $v_2 = 0$ in $\OOO\times
(0,T)$.  That is, $f=0$ in $\OOO$.  By the assumption $f\not\equiv 0$ in
$\OOO$, it turns out that $T - t_* > R(x_0)$ is impossible.
Hence, $T - t_* \le R(x_0)$, that is, $t_* \ge T - R(x_0)$.
Consequently (8.7) implies $\mu(t) = 0$ for $0<t<T-R(x_0)$.
Thus the proof of Proposition B is complete.
\section{Appendix II. Proofs of Lemmas 1 and 2}

We number all the eigenvalues of $A=-\Delta$ with $\mathcal{D}(A)
= H^2(\OOO) \cap H^1_0(\OOO)$ with the multiplicities.
Let $\{ \va_j\}_{j\in \N}$ be eigenfunctions forming an orthonormal
basis in $L^2(\OOO)$: $A\va_j = \sigma_j\va_j$.
Let $T_0>0$ be arbitrary.
\\
{\bf Proof of Lemma 1.}\\
(i) We refer to the example in Evans \cite[Thm. 5 in Ch. 7]{E} to
know that there exists a unique solution $u \in L^2(0,T_0;H^2(\OOO))
\cap H^1(0,T_0; L^2(\OOO))$ to (1.1).
Moreover, similarly to (5.3), we have
$$
u(x,t) = \sumj e^{-\sigma_jt}\left( \int^t_0
e^{\sigma_js}\mu(s) ds \right) (f,\va_j)\va_j(x), \quad x\in \OOO,\,
t>0,                    \eqno{(9.1)}
$$
where the series is convergent in $L^2(0,T_0;H^2(\OOO)) \cap
H^1(0,T_0;L^2(\OOO))$.  Since $f \in C^{\infty}_0(\OOO)
\subset \mathcal{D}(A^{\ell})$ with any $\ell \in \N$ and
$\sigma_j \sim \rho_oj^{\frac{2}{d}}$ as $j \to \infty$ (e.g.,
\cite{Ag}, \cite{CH}), we see (5.6) and
$$
\vert (f,\va_j)\vert \le \frac{C}{j^{\frac{2\ell}{d}}}\Vert A^{\ell}f\Vert
_{L^2(\OOO)}, \quad \ell \in \N.                     \eqno{(9.2)}
$$
Moreover, with $\ell \in \N$ satisfying $\ell_0 > \frac{d}{4}$, the
Sobolev embedding yields
$$
\Vert \va_j\Vert_{C(\ooo{\OOO})} \le C\Vert \va_j\Vert_{H^{2\ell_0}(\OOO)}
\le C\Vert A^{\ell_0}\va_j\Vert_{L^2(\OOO)}
= C\sigma_j^{\ell_0} \le Cj^{\frac{2\ell_0}{d}}.
$$
Therefore, for sufficiently large $\ell \in \N$, we see that the series
in (9.1) is convergent in $L^{\infty}(0,T_0;H^2(\OOO))$, and so
$u \in C([0, T_0];H^2(\OOO))$.  Thus (1.3) is verified.

Since $\ppp_tu = -Au + \mu(t)f$ in $\OOO\times (0,T)$, in terms of
(9.2), we can similarly prove that $\ppp_tu\in C([0, T_0];L^2(\OOO))$.
Thus the proof of Lemma 1 (i) is complete.
\\
(ii) Setting $v = \ppp_tu$, we have
$$
\left\{ \begin{array}{rl}
& \ppp_t^2v = \Delta v + \mu'(t)f(x), \quad x\in \OOO, \, 0<t<T_0,\\
& v(x,0) = 0, \quad \ppp_tv(x,0) = \mu(0)f(x), \quad x\in \OOO, \\
& v(x,t) = 0, \quad x\in \ppp\OOO, \, 0<t<T_0.
\end{array}\right.
$$
Since $\mu' \in L^2(0,T_0)$ and $f \in L^2(\OOO)$, we see that
$v \in C([0,T_0];H^1_0(\OOO)) \cap C^1([0,T_0];L^2(\OOO))$
(e.g., \cite{LM}, Theorem 8.2 (p.275)).
Therefore,
$$
u \in C^1([0, T_0];H^1_0(\OOO)) \cap C^2([0, T_0]; L^2(\OOO)).
$$
Moreover, $\Delta u(\cdot,t) = \ppp_t^2u(\cdot,t) - \mu(t)f$ in $\OOO$ and
$u(\cdot,t) = 0$ on $\ppp\OOO$ for each $t \in [0, T_0]$, and the elliptic
regularity yields $u \in C([0, T_0]; H^2(\OOO))$.
Since
$$
u(x,t) = \sumj \left( \int^t_0 \frac{\sin (t-s)\sqrt{\la_j}}
{\sqrt{\la_j}} \mu(s) ds \right) (f, \va_j)\va_j(x), \quad x\in \OOO,\, t>0
$$
(e.g., \cite{Ko}), in terms of (9.2) we can see that the series
is convergent on $\ooo{\OOO} \times [0, T_0]$, so that
$u \in C(\ooo{\OOO} \times [0, T_0])$ follows.

Finally $\ppp_{\nu}v \in L^2(0, T_0; L^2(\ppp\OOO))$ is
seen by e.g., \cite{Ko}, and so
$\ppp_{\nu}u \in H^1(0, T_0; L^2(\ppp\OOO))$.  Thus the proof of
Lemma 1 is complete.    $\blacksquare$
\\
\vspace{0.2cm}
{\bf Proof of Lemma 2.}\\
The proof is similarly done to Lemma 1 for example by means of
Sections 7.1 and 7.2 of Chapter 7 in Evans \cite{E},
Theorem 8.2 (p.275) in Lions and Magenes \cite{LM}, Komornik \cite{Ko}
and Pazy \cite{Pa}, and we do not repeat the details.

%%%%%%%%%%%%%%%%%%%%%%%%%%%%%%%%%%%%%%%%%%%%%%%%%
%%%%%%%%%%%%%%%%%%%%%%%%%%%%%%%%%%%%%%%%%%%%%%%%%
\section*{Acknowledgments}
The first author was supported by NSFC (no. 11971104, 11971121).
The second author was supported by NSFC (no.11925104) and Program of Shanghai
Academic/Technology Research Leader (19XD1420500).
The third author was supported by Grant-in-Aid (A) 20H00117 of
Japan Society for the Promotion of Science and
by The National Natural Science Foundation of China
(no. 11771270, 91730303).
This paper has been supported by the RUDN University
Strategic Academic Leadership Program.

%%%%%%%%%%%%%%%%%%%%%%%%%%%%%%%%%%%%%%%%%%%%%%%%%
%%%%%%%%%%%%%%%%%%%%%%%%%%%%%%%%%%%%%%%%%%%%%%%%%

\end{document}